\documentclass[11pt,thmsa,a4paper,oneside,final,notitlepage,onecolumn]{article}
\usepackage{amssymb,epsf}
\def\Bbb{\mathbb}
\def\frak{\mathfrak}
\usepackage{sw20jart}
\usepackage{epic}
\usepackage{ecltree}

\input{tcilatex}
\begin{document}

\title{All Abelian Quotient C.I.-Singularities Admit\\
Projective Crepant Resolutions in All Dimensions}
\author{Dimitrios I. Dais, Martin Henk and G\"unter M. Ziegler \\
{\small D.I.D.: Mathematisches Institut der Universit\"{a}t Bonn,
Beringstra\ss e 4, D-53115 Bonn, Germany \ }\\
{\small\tt $\qquad$ dais@mpim-bonn.mpg.de}\\
{\small M.H.: Konrad-Zuse-Zentrum (ZIB), Takustra\ss e 7, D-14195
Berlin-Dahlem, Germany}\\
{\small\tt $\qquad$ henk@zib.de}\\
{\small G.M.Z.: TU-Berlin, FB Mathematik, MA~6-1, Stra\ss e des 17.~Juni 136,
D-10623 Berlin, Germany}\\
{\small\tt $\qquad$ ziegler@math.tu-berlin.de}}
\date{}
\maketitle

\begin{abstract}
\noindent 
For Gorenstein quotient spaces $\Bbb{C}^d/G$, a direct
generalization of the classical McKay correspondence in dimensions $d\geq 4$
would primarily demand the existence of projective, crepant
desingularizations. Since this turned out to be not always possible, Reid
asked about special classes of such quotient spaces which would satisfy the
above property. We prove that the underlying spaces of all Gorenstein
abelian quotient singularities, which are embeddable as complete
intersections of hypersurfaces in an affine space, have torus-equivariant
projective crepant resolutions in all dimensions. We use techniques from
toric and discrete geometry.
\end{abstract}

\section{Introduction}

\noindent Up to isomorphism, the underlying spaces $\Bbb{C}^d/G$, $d\geq 2$,
of Gorenstein quotient singularities can be always realized by finite
subgroups $G$ of SL$\left( d,\Bbb{C}\right) $ acting linearly on $\Bbb{C}^d$. For 
$d=2$, these Gorenstein quotient spaces are embeddable in $\Bbb{C}^3$
as hypersurfaces of $A$-$D$-$E$ type. In 1979, McKay \cite{McK} observed a
remarkable connection between the representation theory of the finite
subgroups of SL$\left( 2,\Bbb{C}\right) $ and the Dynkin diagrams of certain
irreducible root systems. This was the starting-point for
Gonzalez-Sprinberg and Verdier \cite{GSV}, and Kn\"{o}rrer \cite{Kn} to
construct a purely geometric, direct correspondence 
\[
\Psi:\text{Irr}^0\left( G\right) \rightarrow \text{Exc}\left( f\right) 
\]
``of McKay-type'' between the set Irr$^0\left( G\right) $ of non-trivial
irreducible representations of $G$, and the set Exc$\left( f\right) $ of
exceptional prime divisors of the minimal desingularization $f:X\rightarrow 
\Bbb{C}^2/G$, (or, equivalently, between Irr$^0\left( G\right) $ and the
members of the natural basis of the cohomology ring $H^{*}\left( X,\Bbb{Q}\right) 
$). The bijection $\Psi $ induces, in fact, an isomorphism between
the graph of Irr$^0\left( G\right) $ and the dual (resolution-) graph w.r.t.\ 
$f$, i.e., the product of the images of two distinct elements of Irr$^0\left(
G\right) $ under $\Psi $ is mapped onto the exceptional prime divisor
corresponding to the ``right'' graph-vertex. (For various applications of
the geometry of Kleinian singularities, including
also ``quiver-theoretic'' methods, we refer to Slodowy \cite{Slodowy}.
Recently, Ito and Nakamura \cite{IN} gave another interpretation 
of the above correspondence by means of Hilbert schemes
of $G$-orbits.)
\medskip \newline
There are several difficulties (not only of group-theoretic nature) to
generalize (even partially) McKay-type correspondences in higher dimensions.
Already for $d=3$, many things change drastically. For instance, the minimal
embedding dimension of $\Bbb{C}^3/G$ is in many cases very high, and its
singular locus is only rarely a singleton. In addition, \textit{crepant}
desingularizations (i.e., the high-dimensional analogues of the above $f$,
cf.\ \cite{Reid1,Reid3}) are unique only up to ``isomorphism in codimension $1 
$'', and there are lots of examples of non-projective $f$'s. Nevertheless,
Markushevich \cite{Mark}, Ito \cite{Ito1,Ito2,Ito3} and Roan \cite
{Roan1,Roan2} proved case-by-case the \textit{existence} of crepant
resolutions of $\Bbb{C}^3/G$, for all possible finite subgroups $G$ of
SL$\left( 3,\Bbb{C}\right) $, by making use of Blichfeldt's classification
results, and Ito \& Reid \cite{Ito-Reid} established a canonical one-to-one
correspondence between the conjugacy classes of the Tate twist $G\left(
-1\right) $ with age $=1$ and the crepant discrete valuations on $\Bbb{C}^d/G $ 
for any $d\geq 2$. It is worth mentioning that the existence of 
\textit{terminal} Gorenstein quotient singularities for $d\geq 4$ (cf.\ \cite
{MMM,Mo-Ste}) means automatically that \textit{not all }Gorenstein quotient
singularities can be resolved by crepant birational morphisms. On the other
hand, Batyrev \cite{Bat} and Kontsevich \cite{Kon} recently announced proofs
of the invariance of \textit{all} cohomology dimensions of the overlying
spaces of all ``full'' crepant desingularizations for arbitrary $d\geq 2$.
(For the case of abelian acting groups, see \cite[5.4]{BD}). For these
reasons, to proceed the investigation of these quotient spaces in dimensions 
$d\geq 4$ one has \textit{either} to work with at most \textit{partial }crepant 
morphisms (and lose the above quoted cohomology dimension
invariance) \textit{or} to determine \textit{special classes} of quotients
for which this classical type of McKay-correspondence makes sense. As there
are only very few known examples to follow the latter approach to the
problem (cf.\ Hirzebruch-H\"{o}fer \cite[p.~257]{HH} 
and Roan \cite[\S~5]{Roan1}),
Reid asked about conditions on $G$'s which would guarantee the
existence of the required ``full'' crepant (preferably projective)
resolutions of the corresponding quotient spaces (see \cite{Reid4}, 
\cite[\S~4.5]{Ito-Reid}). We believe that one significant class of Gorenstein
quotient spaces which enjoys this property is that of
complete intersections (``c.i.'s'').

\begin{conjecture}
\label{Con}For all finite subgroups $G$ of \emph{SL}$\left( d,\Bbb{C}\right) 
$, for which $\Bbb{C}^d/G$ is minimally embeddable as complete intersection 
\emph{(}in an affine space $\Bbb{C}^r,r\geq d+1$\emph{)}, the quotient space 
$\Bbb{C}^d/G$ admits crepant, projective desingularizations for \emph{all }$d\geq 
2$.
\end{conjecture}

\noindent 
In this paper, we prove the following:

\begin{theorem}
\label{Main}Conjecture \emph{\ref{Con} }is true for all \emph{abelian}
finite groups $G\subset $ \emph{SL}$\left( d,\Bbb{C}\right) $ \emph{(}for which  
$\Bbb{C}^d/G$ is a \emph{``c.i.'').}
\end{theorem}

\noindent For other, non-c.i. abelian quotient spaces having such
desingularizations, we refer to \cite{DHZ}. Our proof makes use of toric
geometry and Watanabe's classification result \cite{Wata}, and is motivated
by the so-called ``nice polyhedral subdivisions'' of Knudsen \& Mumford 
\cite[Ch.~III]{KKMS}. More precisely, the paper is organized as follows: In
\S\ 2 we introduce those notions of the theory of toric varieties which will
enable a convenient combinatorial characterization of the inner structure of 
$\Bbb{C}^d/G$. Useful tools and arguments coming from convex geometry and
the theory of triangulations of polyhedral complexes are presented in \S\ 3.
In section 4 we explain why the existence of projective, crepant resolutions
can be deduced from the existence of b.c.b.-triangulations of the junior
simplex. The method of how one passes from ``special data'' (parametrizing
the defining groups of c.i.-quotient spaces) to ``Watanabe forests'' and to
their lattice-geometric realization via ``Watanabe simplices'' is described
explicitly in \S\ 5. Section 6 provides the inductive, constructive
procedure leading to the desired ``nice'' b.c.b.-triangulations of all
Watanabe simplices. Finally, in \S\ 7 we give concrete formulae for the
computation of the dimensions of the cohomology groups of the spaces
desingularizing $\Bbb{C}^d/G$'s (which are ``c.i'') by any crepant
morphism.\medskip \newline
Concerning the general terminology, we always work with normal complex
varieties, i.e., with normal, integral, separated schemes defined over $\Bbb{C}$. 
A \textit{partial desingularization} of such an $X$ is a proper
birational morphism $f:Y\rightarrow X$ on $X$, in which $Y$ is assumed to be
normal. $f:Y\rightarrow X$ is called a \textit{full} (or 
\textit{overall})\textit{\ desingularization of} $X$ if, in addition, we have 
Sing$\left( Y\right) =\varnothing $, i.e., if the overlying 
space $Y$ is smooth. For brevity's sake, the word 
\textit{desingularization}\emph{\ 
}of $X$ (or the phrase \textit{resolution of singularities of} $X$)
is sometimes used as synonymous with the phrase ``full (or overall)
desingularization of $X$''. When we refer to partial desingularizations, we
mention it explicitly. (By the word \textit{singularity} we intimate either
a singular point or the germ of a singular point, but the meaning will be in
each case clear from the context.)\medskip \newline
A partial desingularization $f:Y\rightarrow X$ of a $\Bbb{Q}$-Gorenstein
complex variety $X$ with global index $j$ is called\ \textit{non-discrepant} or 
simply \textit{crepant} if $\omega _X^{\left[ j\right] }\cong
f_{*}\left( \omega _Y^{\otimes j}\right) $, or, in other words, if the (up
to rational equivalence uniquely determined) difference $jK_Y-f^{*}\left(
jK_X\right) $ contains exceptional prime divisors which have vanishing
multiplicities. ($\omega _X,K_X$ and $\omega _Y,K_Y$ denote the dualizing
sheaves and the canonical divisors of $X$ and $Y$ respectively).
Furthermore, $f:Y\rightarrow X$ is \textit{projective} if $Y$ admits an $f$-ample 
Cartier divisor.

\section{A brief toric glossary\label{TOR-GL}}

\noindent We recall some basic facts from the theory of toric varieties and
fix the notation which will be used in the sequel. For details the reader is
referred to the standard textbooks 
\cite{Ewald} \cite{Fulton} \cite{KKMS} \cite{Oda}.
\medskip

\noindent \textsf{(a)} 
For a set $A$ of vectors of $\Bbb{R}^d$, 
the \textit{linear hull}, the \textit{positive hull}, 
the \textit{affine hull} and \textit{the convex hull} of $A$ is
\begin{eqnarray*}
\mathrm{lin}\left( A\right) &=&\left\{ \sum_{i=1}^k \mu_i \mathbf{x}_i\ \Big|\  
\mathbf{x}_i\in A,\, \mu_i\in \Bbb{R},\, k\in\Bbb{N} \right\} , \\
 \\
\mathrm{pos}\left( A\right) &=&\left\{ \sum_{i=1}^k \mu_i \mathbf{x}_i\ \Big| \  
\mathbf{x}_i\in A,\, \mu_i\in \Bbb{R},\,  \mu_i\ge0,\,k\in\Bbb{N} \right\} , \\
 \\
\mathrm{aff}\left( A\right) &=&\left\{ \sum_{i=1}^k \mu_i \mathbf{x}_i\ \Big|\   
\mathbf{x}_i\in A,\, \mu_i\in \Bbb{R},\, \sum_{i=1}^k \mu_i=1,\,k\in\Bbb{N} 
\right\} ,\ \ \text{and} \\
  \\
\mathrm{conv}\left( A\right) &=&\left\{ \sum_{i=1}^k \mu_i \mathbf{x}_i\ \Big| \  
\mathbf{x}_i\in A,\, \mu_i\in \Bbb{R}_{\geq 0},\, \sum_{i=1}^k 
\mu_i=1,\,k\in\Bbb{N} \right\} ,\ \
\end{eqnarray*}
respectively. Moreover, we define the \textit{integral affine hull} of $A$ as
\[
\text{aff}_{\Bbb{Z}}\left( A\right)  = \left\{ \sum_{i=1}^k \mu_i \mathbf{x}_i\ | 
\  \mathbf{x}_i\in A,\, \mu_i\in \Bbb{Z},\, \sum_{i=1}^k \mu_i=1,\,k\in\Bbb{N} 
\right\}.
\]
For a set $A\subset \Bbb{R}^d$ the  dimension of its affine hull is denoted
by $\dim(A)$.
\medskip
\par\noindent
\textsf{(b) }Let $N\cong \Bbb{Z}^d$ be a free $\Bbb{Z}$-module of rank $d\geq 1$. 
$N$ can be regarded as a \textit{lattice }in $N_{\Bbb{R}}:=N\otimes 
_{\Bbb{Z}}\Bbb{R}\cong \Bbb{R}^d$. (We shall represent the
elements of $N$ by column vectors). If $\left\{ n_1,\ldots ,n_d\right\} $ is
a $\Bbb{Z}$-basis of $N$, then
\[
\text{det}\left( N\right) :=\left| \text{det}\left( n_1,\ldots ,n_d\right)
\right|
\]
is the \textit{lattice determinant}. An $n\in N$ is \textit{primitive}
if conv$\left( \left\{ \mathbf{0},n\right\} \right) \cap N$ contains no
other lattice points except $\mathbf{0}$ and $n$.\smallskip

Let $N\cong \Bbb{Z}^d$ be a lattice as above, $M:=$ Hom$_{\Bbb{Z}}\left( 
N,\Bbb{Z}\right) $ its dual, $N_{\Bbb{R}},M_{\Bbb{R}}$ their real scalar
extensions, and $\left\langle .,.\right\rangle :N_{\Bbb{R}}\times 
M_{\Bbb{R}}\rightarrow \Bbb{R}$ the natural $\Bbb{R}$-bilinear pairing. A subset 
$\sigma $ of $N_{\Bbb{R}}$ is a \textit{strongly convex polyhedral cone}
(\textit{s.c.p.~cone}, for short), if there exist $n_1,\ldots ,n_k\in 
N_{\Bbb{R}}$, such that $\sigma=\mathrm{pos}(\{n_1,\ldots ,n_k\})$,
and $\sigma \cap \left( -\sigma \right) =\left\{ \mathbf{0}\right\} $. Its
\textit{relative interior} int$\left( \sigma \right)$ is the usual topological
interior of it, considered as a subset of lin$\left( \sigma \right) $. The
\textit{dual cone} of an s.c.p.\ cone $\sigma $ is defined by
\[
\sigma ^{\vee }:=\left\{ \mathbf{x}\in M_{\Bbb{R}}\ \left| \ \left\langle
\mathbf{x},\mathbf{y}\right\rangle \geq 0,\ \mbox{ for all } 
\mathbf{y}\in \sigma \right. \right\}
\]
and satisfies: $\sigma ^{\vee }+\left( -\sigma ^{\vee }\right) =M_{\Bbb{R}}$
and dim$\left( \sigma ^{\vee }\right) =d$. \ A subset $\tau $ of a s.c.p.
cone $\sigma $ is called a \textit{face} of $\sigma $ (notation: $\tau \prec
\sigma $), if $\tau =\left\{ \mathbf{y}\in \sigma \ \left| \ \left\langle
m_0,\mathbf{y}\right\rangle =0\right. \right\} $, for some $m_0\in \sigma
^{\vee }$. A s.c.p.\ cone $\sigma$ is
\textit{simplicial} (resp.\ \textit{rational}) if 
$\sigma=\mathrm{pos}(\{n_1,\ldots,n_k\})$, where the vectors 
$n_1,\ldots ,n_k$ are $\Bbb{R}$-linearly independent 
(resp.~if $n_1,\ldots,n_k\in N_{\Bbb{Q}}$). 
If $\sigma \subset N_{\Bbb{R}}$ is a rational s.c.p.\ 
cone, then $\sigma $ is ``pointed'' and the subsemigroup $\sigma \cap N$ of $N$ is 
a monoid having the origin $\mathbf{0}$ as its neutral element. Using
its dual $M\cap \sigma ^{\vee }$, one constructs a finitely generated,
normal $\Bbb{C}$-subalgebra $\Bbb{C}\left[ M\cap \sigma ^{\vee }\right] $ of
$\Bbb{C}\left[ M\right] $ and a $d$-dimensional affine complex variety
\[
U_\sigma :=\text{Max-Spec}\left( \Bbb{C}\left[ M\cap \sigma ^{\vee }\right]
\right) \ .
\]
\noindent \textsf{(c) }For $N\cong \Bbb{Z}^d$ we define an algebraic 
torus\textit{\ }$T_N\cong \left( \Bbb{C}^{*}\right) ^d$ by $T_N:=$ 
Hom$_{\Bbb{Z}}\left( M,\Bbb{C}^{*}\right) =N\otimes _{\Bbb{Z}}\Bbb{C}^{*}$. Every 
$m\in M$
assigns a character $\mathbf{e}\left( m\right) :T_N\rightarrow \Bbb{C}^{*}$
with $\mathbf{e}\left( m\right) \left( t\right) :=t\left( m\right) $, for
all $t\in T_N$. We have
\[
\mathbf{e}\left( m+m'\right) =\mathbf{e}\left( m\right) \cdot
\mathbf{e}\left( m'\right), \ \mbox{ for }m,m'\in M,
\ \ \ \mbox{ and }\mathbf{e}\left( \mathbf{0}\right) =1.
\]
Moreover, for each $n\in N$ we define an $1$-parameter subgroup
of $T_N$ by
\[
\gamma _n:\Bbb{C}^{*}\rightarrow T_N\ \ \ \text{with\ \ \ }\gamma _n\left(
\lambda \right) \left( m\right) :=\lambda ^{\left\langle m,n\right\rangle }\text{, 
\ \ for\ \ \ }\lambda \in \Bbb{C}^{*},\ m\in M,\ \
\]
($\gamma _{n+n'}=\gamma _n\circ $\ $\gamma _{n'}$, \ \
for\ \ \ $n,n'\in N$). We can therefore identify $M$ with the
character group of $T_N$ and $N$ with the group of $1$-parameter subgroups
of $T_N$. Furthermore, $U_\sigma $ (as above) can be identified with the
set of semigroup homomorphisms:
\[
U_\sigma =\left\{ u:M\cap \sigma ^{\vee }\rightarrow \Bbb{C}^*\ \left| \
u\left( \mathbf{0}\right) =1,\ u\left( m+m'\right) =u\left(
m\right) \cdot u\left( m'\right) ,\ \mbox{ for all }
m,m'\in M\cap
\sigma ^{\vee }\right. \right\},
\]
where $\mathbf{e}\left( m\right) \left( u\right):=u\left( m\right)$,
for all $m\in M\cap \sigma ^{\vee }$ and for all $u\in U_\sigma$.
\medskip

\noindent \textsf{(d) }A \textit{fan} $\Delta$ w.r.t.\ $N\cong \Bbb{Z}^d$ is
a finite collection of rational s.c.p.\ cones in $N_{\Bbb{R}}$ 
such that:
\smallskip \newline
(i) any face $\tau $ of $\sigma \in \Delta $ belongs to $\Delta $,
and\smallskip \newline
(ii) for $\sigma _1,\sigma _2\in \Delta $, the intersection $\sigma _1\cap
\sigma _2$ is a face of both $\sigma _1$ and $\sigma _2.\smallskip $\newline
The union $\left| \Delta \right|:=\bigcup \left\{ \sigma \ \left| \ \sigma \in
\Delta \right. \right\} $ is called the \textit{support} of $\Delta $.
Furthermore, we define
\[
\Delta \left( i\right):=\left\{ \sigma \in \Delta \ \left| \ \text{dim}\left( 
\sigma \right) =i\right. \right\},\ \text{for }0\leq i\leq d.
\]
If $\varrho \in \Delta \left( 1\right) $, then there exists a unique
primitive vector $n\left( \varrho \right) \in N\cap \varrho $ with $\varrho 
=\Bbb{R}_{\geq 0}\ n\left( \varrho \right) $ and each cone $\sigma \in \Delta
$ can be therefore written as
\[
\sigma =\sum\Sb \varrho \in \Delta \left( 1\right)  \\ 
\varrho \prec \sigma \endSb \ \Bbb{R}_{\geq 0}\ n\left( \varrho \right).
\]
The set Sk$^1\left( \sigma \right):=\left\{ n\left( \varrho \right) \
\left| \ \varrho \in \Delta \left( 1\right) ,\varrho \prec \sigma \right.
\right\} $ is called the\textit{\ }(\textit{pure}) \textit{first skeleton}
of $\sigma $. The \textit{toric variety X}$\left( N,\Delta \right) $
associated to a fan\textit{\ }$\Delta $ w.r.t.\ the lattice\textit{\ }$N$ is
by definition the identification space
\[
X\left( N,\Delta \right):= \Big( \bigcup_{\sigma \in \Delta }\
U_\sigma \Big) \ \Big/\ \sim
\]
with $U_{\sigma _1}\ni u_1\sim u_2\in U_{\sigma _2}$ if and only if there is
a $\tau \in \Delta ,$ such that $\tau \prec \sigma _1\cap \sigma _2$ and $u_1=u_2$ 
within $U_\tau $. As a complex variety, $X\left( N,\Delta \right) $
turns out to be irreducible, normal, Cohen-Macaulay and to have at most
rational singularities (cf.~\cite[p.~76]{Fulton} and 
\cite[Thm.~1.4, p.~7 and Cor.~3.9, p.~125]{Oda}). 
$X\left( N,\Delta \right) $ admits a canonical $T_N$-action which extends the 
group multiplication of $T_N=U_{\left\{
\mathbf{0}\right\} }$:
\begin{equation}
T_N\times X\left( N,\Delta \right) \ni \left( t,u\right) \longmapsto t\cdot
u\in X\left( N,\Delta \right)   \label{torus action}
\end{equation}
where, for $u\in U_\sigma $, $\left( t\cdot u\right) \left( m\right)
:=t\left( m\right) \cdot u\left( m\right)$, for all $m\in M\cap \sigma
^{\vee }$. The orbits w.r.t.\ the action (\ref{torus action}) are
parametrized by the set of all the cones belonging to $\Delta $. For a $\tau
\in \Delta $, we denote by orb$\left( \tau \right) $ (resp.\ by $V\left( \tau
\right) $) the orbit (resp.\ the closure of the orbit) which is associated to
$\tau $. The spaces orb$\left( \tau \right) $ and $V\left( \tau \right) $
have the following properties (cf.\ \cite[pp.~52-55]{Fulton}, 
\cite[\S~1.3]{Oda}):
\smallskip 
\newline
(i) $V\left( \tau \right) =\coprod \left\{
\text{orb}\left( \sigma \right) \ \left| \ \sigma \in \Delta ,\tau \prec
\sigma \right. \right\} $ and orb$\left( \tau \right) =V\left( \tau \right)
\smallsetminus \bigcup \ \left\{ V\left( \sigma \right) \ \left| \ \tau
\precneqq \sigma \right. \right\} .\smallskip $\newline
(ii) $V\left( \tau \right) =X\left( N\left( \tau
\right) ,\text{ Star}\left( \tau \right) \right) $ is itself a toric variety
w.r.t.
\[
N\left( \tau \right):=N/N_\tau,\ \ \ N_\tau:=N\cap \text{lin}\left(
\Bbb{\tau }\right),\ \ \ \text{Star}\left( \tau \right):=\left\{
\overline{\sigma }\ \left| \ \sigma \in \Delta ,\ \tau \prec \sigma \right.
\right\},
\]
where $\overline{\sigma }\ :=\left( \sigma +\left( N_\tau \right) 
_{\Bbb{R}}\right) \ /\ \left( N_\tau \right) _{\Bbb{R}}$ denotes the image of 
$\sigma
$ in $N\left( \tau \right) _{\Bbb{R}}=N_{\Bbb{R}}\ /\ \left( N_\tau \right)
_{\Bbb{R}}$ .\smallskip \newline
(iii) For $\tau \in \Delta $, $V\left( \tau \right) $ has a natural affine
open covering $\left\{ U_\sigma \left( \tau \right) \ \left| \ \tau \prec
\sigma \right. \right\} $ consisting of ``intermediate'' spaces
\[
U_\tau \left( \tau \right) =\text{ orb}\left( \tau \right) \hookrightarrow
U_\sigma \left( \tau \right) \hookrightarrow U_\sigma
\]
being defined by: $U_\sigma \left( \tau \right) :=$ Max-Spec$\left( \Bbb{C}\left[ 
\overline{\sigma }^{\vee }\cap \ M\left( \tau \right) \right] \right)
$, with $M\left( \tau \right) $ denoting the dual of $N\left( \tau \right) $. It 
should be pointed out, that every $T_N$-invariant subvariety of $U_\sigma $ has 
the form $U_\sigma \left( \tau \right) $ and that dim$\left(
U_\sigma \left( \tau \right) \right) =$ dim$\left( \sigma \right) -$ dim$\left( 
\tau \right) .$
\medskip
\par\noindent
\textsf{(e) }Let $X\left( N,\Delta \right) $ be a $d$-dimensional (not
necessarily compact) simplicial toric variety w.r.t.\ $N\cong{\Bbb Z}^d$. A function $\psi 
:\left| \Delta \right|
\rightarrow \Bbb{R}$ is a  (\textit{rational}) 
$\Delta$-\textit{support function} if $\psi \left( N_{\Bbb{Q}}\cap \left| \Delta \right|\right)
\subset {\Bbb{Q}}$ and $\psi| _\sigma$ is linear for all
$\sigma \in \Delta $. This means that
for all $\sigma\in\Delta$ there has to be some $m_\sigma\in M_{\Bbb{Q}}$
such that
\[
\psi ( \mathbf{x}) =\left\langle m_\sigma ,\mathbf{x}\right\rangle
\ \ \mbox{for all } \mathbf{x}\in \sigma
\qquad\mbox{and}\qquad
\left\langle m_\sigma ,\mathbf{x}\right\rangle = 
\left\langle m_\tau   ,\mathbf{x}\right\rangle
\mbox{\ whenever\ } \tau\prec \sigma\mbox{\ and\ }\mathbf{x}\in\tau.
\]
We denote by SF$_{\Bbb{Q}}\left( N,\Delta \right) $ the additive group of all 
rational $\Delta$-support functions. 
A $\psi \in $ SF$_{\Bbb{Q}}\left( N,\Delta \right) $ is 
\textit{strictly upper convex} if for
every maximal $\sigma \in \Delta $  the $m_\sigma$ can be chosen such that 
$\psi( \mathbf{x} ) \leq \left\langle m_\sigma ,\mathbf{x}\right\rangle $,
with equality if and only if $\mathbf{x}\in \sigma$. Let
\[
\text{SUCSF}_{\Bbb{Q}}\left( N,\Delta \right) :=\left\{ \psi \in \text{SF}_{\Bbb{Q}}\left(
N,\Delta \right) \ \left| \ \psi \ \text{ strictly upper convex}\right.
\right\}.
\]
The group $T_N$-CDiv$_{\Bbb{Q}}\left( X\left( N,\Delta \right) \right) $ of $T_N$-invariant 
${\Bbb{Q}}$-Cartier divisors on $X\left( N,\Delta \right) $ has the $\Bbb{Q}$-basis $\left\{ 
\Bbb{Q\ }V\left( \varrho \right) \ \left| \  \varrho \in \Delta \left( 1\right)
\right. \right\} $  (cf.\ \cite[p.~68-69]{Oda}).

\begin{theorem}
\label{SUCF}The relationship between the rational $\Delta $-support functions and the
$T_N$-invar\-iant $\Bbb{Q}$-Cartier divisors on $X\left( N,\Delta \right) $ is given by
the following bijections:
\begin{equation}
\begin{array}{cccccc}
\psi  & \in  & \text{\emph{SF}}_{\Bbb{Q}}\left( N,\Delta \right)  &  & \supset  &
\text{\emph{SUCSF}}_{\Bbb{Q}}\left( N,\Delta \right)  \\
\updownarrow  &  & \updownarrow \ \text{\emph{1:1}} &  &  & \updownarrow \
\text{\emph{1:1}} \\
D_\psi  & \in  & T_N\text{-\emph{CDiv}}_{\Bbb{Q}}\left( X\left( N,\Delta \right)
\right)  & =\dbigoplus_{\varrho \in \Delta \left( 1\right) }\ \Bbb{Q\ }V\left( 
\varrho \right)  & \supset  & \left\{
\begin{array}{c}
\ \text{\emph{ample $\Bbb{Q}$-Cartier }} \\
\text{\emph{\  divisors on \ }}X\left( N,\Delta \right) \
\end{array}
\right\}
\end{array}
\label{div-sup}
\end{equation}
with
\[
D_\psi:=-\sum_{\varrho \in \Delta \left( 1\right) }\ \psi \left( n\left(
\varrho \right) \right) \ V\left( \varrho \right).
\]
\end{theorem}

\noindent \textit{Proof. }
For the first bijection see \cite[Prop.\ 2.1, p.~68-69]{Oda}. 
The proof of the second bijection of (\ref{div-sup}), in the case
in which $X\left( N,\Delta \right) $ is compact, is given in 
\cite[Cor.~2.14, p.~83]{Oda}. 
The general case is treated in \cite[Ch.~I, Thm.~9,
p.~42, and Thm. 13, p.~48]{KKMS}. $_{\Box }$

\begin{corollary}
\label{quasiproj}Let $X\left( N,\Delta \right) $ be a simplicial toric variety \emph{(}resp.\ 
a simplicial compact toric variety\emph{).} Then $X\left( N,\Delta \right) $ is
quasiprojective \emph{(}resp.\ projective\emph{)} if and only if 
\emph{SUCSF}$\left( N,\Delta \right) \neq \varnothing .$
\end{corollary}

\noindent \textsf{(f)}
A \textit{map of fans}
$\varpi :( N',\Delta ')\rightarrow (N,\Delta) $
is a $\Bbb{Z}$-linear homomorphism $\varpi:N'\rightarrow N$
whose extension $\varpi_{\Bbb{R}}:N'_{\Bbb{R}}\rightarrow N_{\Bbb{R}}$
satisfies the property
\[
\text{for all }\sigma'\in\Delta'
\text{ there is some }\sigma\in\Delta
\text{ such that }\varpi_{\Bbb{R}}(\sigma')\subseteq\sigma.
\]
Every map of fans 
$\varpi :( N',\Delta ')\rightarrow (N,\Delta) $ induces an holomorphic map
$$
\varpi _{*}:X\left( N',\Delta
'\right) \rightarrow X\left( N,\Delta \right) 
$$
which is equivariant
w.r.t.\ the actions of $T_{N'}$ and $T_N$ on the toric varieties $X\left( 
N',\Delta '\right) $, $X\left( N,\Delta \right) $.

\begin{theorem}
\label{prop-bir}If $\varpi :\left( N',\Delta '\right)
\rightarrow \left( N,\Delta \right) $ is a map of fans, $\varpi _{*}$ is
\emph{proper} if and only if $\varpi ^{-1}\left( \left| \Delta \right|
\right) =\left| \Delta '\right| .$ In particular, if $N=N'$
and $\Delta '$ is a refinement of $\Delta $, i.e., if each cone of $\Delta $ is a 
union of cones of $\Delta '$, then \emph{id}$_{*}:X\left( N,\Delta '\right) 
\rightarrow X\left( N,\Delta \right)
$ is \emph{proper }and \emph{birational.}
\end{theorem}

\noindent \textit{Proof.}
See \cite[Thm.~1.15, p.~20, and Cor.~1.18, p.~23]{Oda}. 
$_{\Box }\bigskip $
\newline
\textsf{(g) }Let $N\cong \Bbb{Z}^d$ be a lattice of rank $d$ and $\sigma
\subset N_{\Bbb{R}}$ a simplicial, rational s.c.p.\ cone of dimension $k\leq d$.
$\sigma $ can be obviously written as $\sigma =\varrho _1+\cdots +\varrho _k$, for 
distinct $1$-dimensional cones $\varrho _1,\ldots ,\varrho _k$. We
denote by
\[
\mathbf{Par}\left( \sigma \right) :=
\Big\{ \mathbf{x}\in ( N_\sigma) _{\Bbb{R}}\ \Big| \ 
\mathbf{x}=\sum_{j=1}^k\ \varepsilon _j\
n( \varrho _j) ,\ \text{with\ \ }0\leq \varepsilon _j<1\mbox{ for }
1\leq j\leq k \Big\}
\]
the \textit{fundamental }(\textit{half-open})\textit{\ parallelotope }which
is associated to $\sigma $. The \textit{multiplicity} mult$\left(
\sigma ,N\right) $ of $\sigma $ with respect to $N$ is defined as
\[
\text{mult}\left( \sigma ,N\right) :=\#\left( \mathbf{Par}\left( \sigma
\right) \cap N_\sigma \right) =\text{Vol}\left( \mathbf{Par}\left( \sigma
\right) ;N_\sigma \right),
\]
where
\[
\text{Vol}\left( \mathbf{Par}\left( \sigma \right) ;N_\sigma \right) 
:=\frac{\text{Vol}\left( \mathbf{Par}\left( \sigma \right) \right) 
}{\text{det}\left( N_\sigma \right) }
\]
is the relative volume of $\mathbf{Par}( \sigma) $ w.r.t.\ $N_\sigma $.
The affine toric variety $U_\sigma $ is smooth if and only if
mult$( \sigma ,N) =1$.

\begin{theorem}
\label{TORRES}Every toric variety $X\left( N,\Delta \right) $ admits a 
$T_N$-equivariant desingularization
\[
f=\emph{id}_{*}:X\left( N,\Delta '\right) \rightarrow X\left(
N,\Delta \right)
\]
by a suitable refinement $\Delta '$ of $\Delta $.
\end{theorem}

\noindent 
In fact, one can always make $\Delta$ simplicial
without introducing additional rays. (This step relies on
Carath\'{e}odory's theorem, cf.\ 
\cite[III 2.6, p.~75]{Ewald}). 
\smallskip 
In a second step, this new simplicial $\Delta $ will be
subdivided further into subcones of strictly smaller multiplicities than
those of the cones of the starting-point. (Since for each $\sigma \in \Delta
$, mult$\left( \sigma ,N\right) $ is a volume, mult$\left( \widehat{\sigma 
},N\right) <$ mult$\left( \sigma ,N\right) $ for every simplicial subcone 
$\widehat{\sigma }$ of $\sigma $). $\Delta '$ is constructed after
finitely many subdivisions of this kind (cf.\ \cite[pp.~31-35]{KKMS}).
\bigskip \newline
\textsf{(h) }For the germ $\left( U_\sigma ,\text{orb}\left( \sigma \right)
\right) $ of an affine $d$-dimensional toric variety w.r.t.\ a \textit{singular 
point }orb$\left( \sigma \right) $, (with dim$\left( \sigma \right)
=d=$ dim$\left( N\right) $), we define the \textit{splitting codimension}
splcod$\left( \text{orb}\left( \sigma \right) ;U_\sigma \right) $ of orb$\left( 
\sigma \right) $ in $U_\sigma $ as
\[
\text{splcod}\left( \text{orb}\left( \sigma \right) ;U_\sigma \right) :=\text{max 
}\left\{ \varkappa \in \left\{ 2,\ldots ,d\right\} \
\left| \
\begin{array}{c}
U_\sigma \cong U_{\sigma '}\times \Bbb{C}^{d-\varkappa },\ \text{dim}\left( \sigma 
'\right) =\varkappa , \\
\text{and\ \ \ Sing}\left( U_{\sigma '}\right) \neq \varnothing
\end{array}
\right. \right\}.
\]
(For orb$\left( \sigma \right) $ regular we can formally define splcod$\left( 
\text{orb}\left( \sigma \right) ;U_\sigma \right) =0$.) If
splcod$\left( \text{orb}\left( \sigma \right) ;U_\sigma \right) = d$, then 
orb$\left( \sigma \right) $\ will be called an \textit{msc-singularity}, i.e., a
singularity having the maximum splitting codimension. An msc-singularity $\left( 
U_\sigma ,\text{orb}\left( \sigma \right) \right) $ 
is \textit{absolutely unbreakable} if $U_\sigma $ cannot be
isomorphic even to the product of (at least two) \textit{singular }affine
toric varieties.\bigskip

\noindent \textsf{(i) }For any abelian finite $G\subset $ GL$\left( 
d,\Bbb{C}\right) $, $d\geq 2$, of order $l\geq 2$, which is \textit{small} (i.e.
which has no pseudoreflections), $\Bbb{C}^d/G$ is singular, with singular
locus Sing$\left( \Bbb{C}^d/G\right)$ containing at least the image $\left[
\mathbf{0}\right] $ of the origin under the canonical quotient-map 
$\Bbb{C}^d\rightarrow \Bbb{C}^d/G$. If we fix a decomposition
\[
G\cong \left( \Bbb{Z}/q_1\Bbb{Z}\right) \times \cdots \times \left( 
\Bbb{Z}/q_\kappa \Bbb{Z}\right)
\]
into cyclic groups, take primitive
$q_i$th roots of unity $\zeta _{q_i}$, and choose
eigencoordinates $z_1,\ldots ,z_d$, then the action of $G$ on $\Bbb{C}^d$ is
determined by
\[
\mathrm{GL}(d,{\Bbb C})\times{\Bbb C}^d \hookleftarrow
\left( \Bbb{Z}/q_i\Bbb{Z}\right) \times \Bbb{C}^d\ni \left( g_i,\left(
z_1,\ldots ,z_d\right) \right) \longmapsto \left( \zeta
_{q_i}^{\alpha _{i,1}}z_1,\ldots ,\zeta _{q_i}^{\alpha _{i,d}}z_d\right) \in
\Bbb{C}^d,\ \ 1\leq i\leq \kappa ,
\]
for generators $g_i=$ diag$\left( \zeta _{q_i}^{\alpha _{i,1}},\ldots ,\zeta
_{q_i}^{\alpha _{i,d}}\right) $, where the weights $\left( \alpha
_{i,1},\ldots ,\alpha _{i,d}\right) $ are unique up to only the usual
conjugacy relations. Using the above toric glossary, we may identify $\Bbb{C}^d/G$ 
with the affine toric variety $U_{\sigma _0}=$ Max-Spec$\left( \Bbb{C}\left[ 
\sigma _0^{\vee }\cap M_G\right] \right) $ being associated to the
positive orthant $\sigma _0=
\mathrm{pos}(\{e_1,\ldots,e_d\})$
with respect to the lattice of weights
\[
N_G=\Bbb{Z}^d+\sum_{i=1}^\kappa \ \Bbb{Z\ }\left( \frac 1{q_i}\left( \alpha
_{i,1},\ldots ,\alpha _{i,d}\right)^\intercal \right) \ \ \ \ 
\text{with \ det}( N_G) =\frac 1l,
\]
where $M_G$ denotes the dual of $N_G$,  which can be viewed as the lattice
parametrizing all $G$-invariant Laurent monomials. Formally, we identify 
$\left[ \mathbf{0}\right] $ with orb$\left( \sigma _0\right) $ and
$U_{\sigma _0}$ with the toric variety $X\left( N_G,\Delta _G\right) $
w.r.t.\ the fan $\Delta _G$ consisting of $\sigma _0$ itself, together
with all its faces. In
these terms, and since $\sigma _0$ is a simplicial s.c.p.\ cone, the singular
locus of $X\left( N_G,\Delta _G\right) $ can be written as the union
\[
\text{Sing}\left( X\left( N_G,\Delta _G\right) \right) =\text{ orb}\left(
\sigma _0\right) \cup \left( \bigcup \left\{ U_{\sigma _0}\left( \tau
\right) \,\left| \,\tau \precneqq \sigma _0,\ 
\mathrm{mult}(\tau, N_G) \geq 2 
 \right. \right\} \right).
\]
For $0\leq i\leq d$, $\Delta _G\left( i\right) $ has\thinspace $\binom
di\,$simplicial cones. In particular, $\Delta _G\left( 0\right) =\left\{
\mathbf{0}\right\} $, and for $1\leq i\leq d$,
\[
\Delta _G\left( i\right) =\left\{
\sigma _0\left( \nu
_1,\ldots ,\nu _i\right) :=
\mathrm{pos}(\{e_{\nu _1},\ldots,e_{\nu _i}\})
\left| \,1\leq \nu _1<\nu _2<\cdots <\nu _i\leq
d\right. \right\} \,.
\]
For any cone $\sigma _0\left( \nu _1,\ldots ,\nu _i\right) \in \Delta
_G\left( i\right) $, $U_{\sigma _0}\left( \sigma _0\left( \nu _1,\ldots ,\nu
_i\right) \right) $ is nothing but
\[
( \{ \mathbf{z}=( z_1,\ldots ,z_d) \in \Bbb{C}^d\,| \,z_{\nu _1}=\cdots =\,z_{\nu 
_i}=0 \} ) /G\,,
\]
i.e., a $\left( d-i\right) $-dimensional coordinate-subspace of $\Bbb{C}^d$
divided by the inherited $G$-action.\smallskip

\begin{proposition}
\label{Gor-prop}For an abelian $G$ \emph{(}as above\emph{)} the following
conditions are equivalent \emph{:\smallskip }\newline
\emph{(i)} $X\left( N_G,\Delta _G\right) =U_{\sigma _0}=\Bbb{C}^d/G$ is
Gorenstein,\smallskip \newline
\emph{(ii) }$G\subset $ \emph{SL}$\left( d,\Bbb{C}\right) $,\smallskip
\newline
$\smallskip $\emph{(iii)} $\sum_{j=1}^d\ \alpha _{i,j}\equiv  0\ 
\text{\emph{mod}}\ q_i  $, for $1\leq i\leq \kappa $,\smallskip \newline
\emph{(iv)} $\left\langle ( 1,1,\ldots ,1) ,n\right\rangle \geq
1$, for all $n\in \sigma _0\cap \left( N_G\smallsetminus \left\{ 
\mathbf{0}\right\} \right) $,\smallskip \newline
\emph{(v) }$\left( U_{\sigma _0},\text{\emph{orb}}\left( \sigma _0\right)
\right) $ is a canonical singularity of index $1$.
\end{proposition}

\noindent \textit{Proof.} (i) $\Leftrightarrow $ (ii) follows from 
\cite{Wat} and (i) $\Leftrightarrow $ (v) from \cite{Reid1,Reid2}. All the other
implications can be checked easily. $_{\Box }\bigskip $

\noindent 
{From} now on, let $X\left( N_G,\Delta _G\right) $ be Gorenstein.
The cone $\sigma _0=\Bbb{R}_{\geq 0}\,\frak{s}_G$ is supported by the
so-called 
\textit{junior lattice simplex }$\frak{s}_G=$ conv$\left( \left\{
e_1,\ldots,e_d\right\} \right) $ (w.r.t.\ $N_G$;
cf.\ \cite{Ito-Reid} \cite{BD}). Note that up to $\mathbf{0}$ there is no
other lattice point of $\sigma _0\cap N_G $ lying ``under'' the affine
hyperplane of $\Bbb{R}^d$ containing $\frak{s}_G$.  Moreover, the
lattice points representing the $l-1$ non-trivial group elements are
exactly those belonging to the intersection of a dilation 
$\lambda\,\frak{s}_G$ of $\frak{s}_G$ with 
$\mathbf{Par}\left( \sigma _0\right)$, for some integer $\lambda$,
$1\le\lambda\le d-1$. For c.i.-$U_{\sigma _0}$'s,
our purpose is to construct crepant, projective (full)
desingularizations 
\[
f:X\left( N_G,\widehat{\Delta _G}\right)
\rightarrow X\left( N_G,\Delta _G\right) =U_{\sigma _0}
\]
which are
$T_{N_G}$-equivariant (with $T_{N_G}=\left( \Bbb{C}^{*}\right) ^d/G$
denoting the algebraic torus embedded in~$U_{\sigma _0}$), and have
overlying (quasiprojective) spaces $X\left( N_G,\widehat{\Delta
_G}\right) $ associated to fans $\widehat{\Delta _G}$ which refine
appropriately $\Delta _G$.

\section{On b.c.b.-triangulations of lattice polytopes}

\noindent In this section we introduce ``b.c.b.-triangulations''
and study their behaviour 
with respect to \textit{joins} and \textit{dilations}. 
(We shall mostly use the standard terminology from the theory of polyhedral
complexes and polytopes, cf.\ \cite{Stanley} \cite{Ziegler}).\bigskip

\noindent \textsf{(a) }By vert$\left( \mathcal{S}\right) $ we denote the set
of vertices of a polyhedral complex $\mathcal{S}$. 
By a \textit{triangulation}~$\mathcal{T}$ of a polyhedral complex 
$\mathcal{S}$ we mean a
geometric simplicial subdivision of $\mathcal{S}$ with vert$\left( 
\mathcal{S}\right) \subset $ vert$\left( \mathcal{T}\right) $. If $\mathcal{T}$ is 
$d$-dimensional, then by $\mathcal{T}( i)$, $0\le i\le d$, we
denote the set of $i$-dimensional simplices of $\mathcal{T}$. A polytope $P$
will be frequently identified with the polyhedral complex consisting of $P$
itself together with all its faces.\bigskip

\noindent \textsf{(b) }A triangulation $\mathcal{T}$ of a polyhedral $d$-complex 
$\mathcal{S}$ in ${\Bbb R}^d$ is \textit{coherent }(or \textit{regular}) if
there exists a strictly upper convex $\mathcal{T}$-support function 
$\psi:\left| \mathcal{T}\right| \rightarrow \Bbb{R}$, i.e., 
a piecewise linear real function
on the underlying space $\left| \mathcal{T}\right| $ for which
\[
\psi ( t\ \mathbf{x}+( 1-t) \mathbf{y})\ \ge\ t\ \psi
( \mathbf{x}) +( 1-t) \psi \left( \mathbf{y}\right)
\text{ \ for all }\mathbf{x},\mathbf{y}\in \left| \mathcal{T}\right|
\text{\ and\ }\emph{\ }t\in \left[ 0,1\right] ,
\]
such that for each maximal simplex $\mathbf{s}$
there is a linear function 
$h_{\mathbf{s}}(\mathbf{x})$
satisfying $\psi(\mathbf{x}) \le h_{\mathbf{s}}(\mathbf{x})$ for all
$\mathbf{x}\in|\mathcal{T}|$,
with equality if and only if $\mathbf{x}\in\mathbf{s}$.
The set of all (real) strictly upper convex $\mathcal{T}$-support
functions is denoted by SUCSF$_{\Bbb{R}}\left( \mathcal{T}\right) $.

\begin{lemma}[Patching Lemma]
\label{PATCH}Let $P\subset \Bbb{R}^d$ be a $d$-polytope, $\mathcal{T}=\left\{ 
\mathbf{s}_i\ \left| \ i\in I\right. \right\} \ $\emph{(}with $I$ a
finite set\emph{) }a coherent triangulation of $P$, and $\mathcal{T}_i=\left\{ 
\mathbf{s}_{i,j}\ \left| \ j\in J_i\right. \right\} \ $\emph{(}$J_i$ finite, for 
all $i\in I$\emph{) }a coherent triangulation of $\mathbf{s}_i$, for all $i\in I$. 
If $\psi _i:\left| \mathcal{T}_i\right| \rightarrow
\Bbb{R}$ denote strictly upper convex $\mathcal{T}_i$-support functions,
such that
\[
\psi _i\left| _{\mathbf{s}_i\cap \mathbf{s}_{i'}}\right. =\psi
_{i'}\left| _{\mathbf{s}_i\cap \mathbf{s}_{i'}}\right.
\]
for all $\left( i,i'\right) \in I\times I$, then
\[
\widetilde{\mathcal{T}}:=\{ \text{\emph{all the simplices }}\mathbf{s}_{i,j}\ 
| \ \ j\in J_i,\ i\in I\}
\]
forms a \emph{coherent triangulation }of the initial polytope $P$.
\end{lemma}

\noindent 
The above $\psi _i$'s can be canonically ``patched together''
to construct an element $\psi$ of {SUCSF}$_{\Bbb{R}}(\widetilde{\mathcal{T}})$;
see \cite[Cor.~1.12, p.~115]{KKMS} or \cite[Lemma\ 3.2.2]{BGT}.\bigskip

\noindent \textsf{(c) }A triangulation $\mathcal{T}$ of a $d$-dimensional
simplex $\mathbf{s}$ (or, more generally, of any pure $d$-dimensional
simplicial complex $\mathcal{S}$) is \textit{balanced} if its graph
can be $( d+1)$-coloured, i.e., if there is a function
\[
\varphi :\text{vert}\left( \mathcal{T}\right) \rightarrow \left\{
0,1,2,\ldots ,d\right\}
\]
such that any two adjacent vertices receive different values
(``colours'') under $\varphi$. If $\mathcal{T}$ is a 
balanced triangulation of a
$d$-polytope, then all facets of $\mathcal{T}$ receive all
the colours, and the colouring function $\varphi$ is unique (up to
permutation of the colours).

\begin{example}
{\rm
The first two $2$-dimensional triangulations $\mathbf{(A)}$ and $\mathbf{(B)}$
of Figure 1 are balanced, whereas $\mathbf{(C)}$ and $\mathbf{(D)}$ are not.}
\end{example}

\begin{center}
\unitlength 0.6truecm
\begin{picture}(18,3.5)
%
\put(0,0){\circle*{0.15}}
\put(1,0){\circle*{0.15}}
\put(2,0){\circle*{0.15}}
\put(3,0){\circle*{0.15}}
\put(1,1){\circle*{0.15}}
\put(2,1){\circle*{0.15}}
\put(3,1){\circle*{0.15}}
\put(2,2){\circle*{0.15}}
\put(3,2){\circle*{0.15}}
\put(3,3){\circle*{0.15}}
\thinlines
\put(0,0){\line(1,0){3}} \put(0,-0.5){\bf{0}}
\put(1,1){\line(1,0){2}} \put(1,-0.5){\bf{1}}
\put(2,2){\line(1,0){1}} \put(2,-0.5){\bf{2}}
\put(3,0){\line(0,1){3}} \put(3.2,-0.5){\bf{0}}
\put(2,0){\line(0,1){2}} \put(0.5,1){\bf{2}}
\put(1,0){\line(0,1){1}} \put(2.1,1){\bf{0}}
\put(0,0){\line(1,1){3}} \put(3.2,1){\bf{1}}
\put(1,0){\line(1,1){2}} \put(1.5,2){\bf{1}}
\put(2,0){\line(1,1){1}} \put(3.2,2){\bf{2}}
\put(3.2,3){\bf{0}}
\put(1,-1.5){$\mathbf{(A)}$}
\put(8,1){
\put(0,-1){\circle*{0.15}}
\put(-1,0){\circle*{0.15}}
\put(0,0){\circle*{0.15}}
\put(1,0){\circle*{0.15}}
\put(0,1){\circle*{0.15}}
\thinlines
\put(-1,0){\line(1,0){2}}
\put(0,-1){\line(0,1){2}}
\put(1,0){\line(-1,1){1}}
\put(0,-1){\line(-1,1){1}}
\put(-1,0){\line(1,1){1}}
\put(0,-1){\line(1,1){1}}
\put(1.2,-0.2){\bf{0}}
\put(-1.5,-.2){\bf{0}}
\put(0.1,0.1){\bf{2}}
\put(-0.1,1.2){\bf{1}}
\put(-0.1,-1.6){\bf{1}}
\put(-0.5,-2.5){$\mathbf{(B)}$}
}
\put(12,1){
\put(-1,-1){\circle*{0.15}}
\put(0,0){\circle*{0.15}}
\put(1,0){\circle*{0.15}}
\put(0,1){\circle*{0.15}}
\thinlines
\put(0,0){\line(1,0){1}}
\put(0,0){\line(0,1){1}}
\put(0,0){\line(-1,-1){1}}
\put(-1,-1){\line(1,2){1}}
\put(-1,-1){\line(2,1){2}}
\put(1,0){\line(-1,1){1}}
\put(-0.5,-2.5){$\mathbf{(C)}$}
}
\put(15,0){
\put(0,0){\circle*{0.15}}
\put(1,0){\circle*{0.15}}
\put(2,0){\circle*{0.15}}
\put(3,0){\circle*{0.15}}
\put(1,1){\circle*{0.15}}
\put(2,1){\circle*{0.15}}
\put(3,1){\circle*{0.15}}
\put(2,2){\circle*{0.15}}
\put(3,2){\circle*{0.15}}
\put(3,3){\circle*{0.15}}
\thinlines
\put(0,0){\line(1,0){3}}
\put(1,1){\line(1,0){2}}
\put(2,2){\line(1,0){1}}
\put(3,0){\line(0,1){3}}
\put(2,0){\line(0,1){2}}
\put(1,0){\line(0,1){1}}
\put(0,0){\line(1,1){3}}
\put(1,0){\line(1,1){1}}
\put(3,1){\line(-1,1){1}}
\put(2,0){\line(1,1){1}}
\put(1,-1.5){$\mathbf{(D)}$}
}
\end{picture}
\vspace{0.8truecm}
\begin{center} {\bf Figure 1} \end{center}
\end{center}

\noindent \textsf{(d) }Let $N$ be a $d$-dimensional lattice. 
A \textit{lattice polytope} (w.r.t.\  $N$) is
a polytope in $N_{\Bbb{R}}\Bbb{\cong R}^d$ all of whose vertices belong to $N$. 
If $\{n_0,n_1,\ldots ,n_k\} $ is a set of affinely independent
lattice points ($k\leq d$), $\mathbf{s}$ the $k$-simplex $\mathbf{s}=$ 
conv$\left( \left\{ n_0,n_1,\ldots ,n_k\right\} \right)$, and 
$N_{\mathbf{s}}:=$ lin$( \{ n_1-n_0,\ldots ,n_k-n_0\} ) \cap N$,
then $\mathbf{s}$ is \textit{basic} if it satisfies any one of
the following equivalent conditions:
\smallskip \newline
(i) $\text{\rm aff}_{\Bbb{Z}}(n_0,\dots,n_k)=N\cap\text{\rm aff}(\mathbf{s})$,
\smallskip\newline
(ii) $\left\{ n_1-n_0,n_2-n_0,\ldots ,n_k-n_0\right\} $ is a $\Bbb{Z}$-basis
of $N_{\mathbf{s}}$,\smallskip \newline
(iii) $\left\{ n_1-n_0,n_2-n_0,\ldots ,n_k-n_0\right\} $ is a part of a 
$\Bbb{Z}$-basis of $N$,\smallskip \newline
(iv) $\mathbf{s\ }$has relative volume\ Vol$\left( 
\mathbf{s};N_{\mathbf{s}}\right) =\dfrac{\text{Vol}\left( \mathbf{s}\right) 
}{\text{det}\left( N_{\mathbf{s}}\right) }=\dfrac 1{k!}\ \ $(w.r.t.\ 
$N_{\mathbf{s}}$).
\smallskip

\noindent
A \textit{lattice triangulation }$\mathcal{T}$ of a lattice polytope $P$ 
(w.r.t.\ $N\cong \Bbb{Z}^d$) is a triangulation of $P$ with
vertices in $N$, and it is called \textit{basic} if all its simplices are
basic. If\/ $\mathcal{T}$ is a basic triangulation of $P$ w.r.t.~$N$, then
in particular we have 
\[
\mathrm{aff}_{\Bbb{Z}}(\mathrm{vert}(\mathcal{T}))\ =\ 
N \cap \mathrm{aff}(P).
\]
Affine maps $\Phi:N_{\Bbb{R}}\rightarrow N_{\Bbb{R}}$
will be called 
\textit{affine integral transformations}
(w.r.t\ $N$) if they satisfy $N=\Phi(N)$,
that is, if they have the form
$\Phi \left(
\mathbf{x}\right) =\mathcal{M}\mathbf{x}+\beta $, with $\mathcal{M}\in $ 
GL$_{\Bbb{Z}}\left( N\right) \cong $ GL$\left( d,\Bbb{Z}\right) $ and $\beta
\in N$ (i.e., if they are composed of $N$-unimodular transformations and 
$N$-integral translations). 


\begin{definition}
\emph{From now on, \textit{b.c.b.-triangulation}
is used as an abbreviation for a basic, coherent, balanced triangulation 
of a lattice polytope.} 
\end{definition}

\begin{lemma}[Preservation of the ``b.c.b.-property'' 
under affine  transformations] \hfill\newline 
Let $\mathcal{T}$ be a lattice triangulation of a lattice polytope $P\subset 
N_{\Bbb{R}}$ \emph{(}w.r.t.\ $N\cong \Bbb{Z}^d$\emph{)} and let 
\mbox{$\Phi :N_{\Bbb{R}}\rightarrow N_{\Bbb{R}}$} 
be a \hbox{\rm (}not necessarily integral\hbox{\rm )} regular
affine transformation. 
If\/\ $\mathcal{T}$ is a b.c.b.-triangulation w.r.t.~$N$, then its image $\Phi 
\left( \mathcal{T}\right) :=\left\{ \Phi \left( \mathbf{s}\right) \ \left| \ 
\mathbf{s}\in
\mathcal{T}\right. \right\} $ under $\Phi $ is again a b.c.b.-triangulation
of the transformed lattice $d$-polytope $\Phi \left( P\right) $ w.r.t.\ the 
lattice $\Phi(N)$.
\label{AFITR}
\end{lemma}


\noindent \textsf{(e) }If $\mathcal{T}_1$ (resp.\ $\mathcal{T}_2$) is a
triangulation of a $d_1$-polytope $P_1$ (resp.\ of a $d_2$-polytope $P_2$) in
$\Bbb{R}^d$ such that 
$
\dim(P_1\cup P_2)=\dim(P_1)+\dim(P_2)+1,
$
that is, such that 
$\mathrm{aff}(P_1)$ and 
$\mathrm{aff}(P_2)$ 
are skew affine subspaces of~$N_{\Bbb{R}}$,
then the \textit{join }$\mathcal{T}_1*\mathcal{T}_2$ of $\mathcal{T}_1$ and 
$\mathcal{T}_2$ is defined~by
$$
\mathcal{T}_1*\mathcal{T}_2=\left\{ \text{conv}\left( \mathbf{s}_1\cup 
\mathbf{s}_2\right) \,\left| \,\mathbf{s}_1\in \mathcal{T}_1,\mathbf{s}_2\in 
\mathcal{T}_2\right. \right\}.
$$
It is a triangulation of the
usual join $P=$ conv$\left( P_1\cup P_2\right) $ of $P_1$ with $P_2$,
which is a polytope of dimension $d_1+d_2+1$ (cf.~\cite[Ex.~7,
p.~136]{Stanley} and \cite[Ex.~9.9, p.~323]{Ziegler}). 


\begin{theorem}[Preservation of the ``b.c.b.-property'' for joins]
\hfill \newline 
Let $\mathcal{T}_1,\mathcal{T}_2$ be triangulations of polytopes
$P_1$, $P_2\subset N_{\Bbb R}\cong\Bbb{R}^d$ of dimensions $d_1$ resp.\ $d_2$
in skew affine subspaces of $N_{\Bbb{R}}$, i.e., satisfying 
$
\dim(P_1\cup P_2)=d_1+d_2+1.
$

If \ $\mathcal{T}_1$, $\mathcal{T}_2$ are b.c.b.~triangulations 
of $P_1$,  $P_2$ w.r.t.~$N$,
then $\mathcal{T}_1*\mathcal{T}_2$ is a
b.c.b.-triangulation of the $(d_1+d_2+1)$-polytope 
$P=\mbox{\rm conv}(P_1\cup P_2)$ w.r.t.~$N$
if and only if 
\[
\mathrm{aff}_{\Bbb{Z}}\big(\,
\mathrm{vert}(\mathcal{T}_1)\cup 
\mathrm{vert}(\mathcal{T}_2)\,\big) \ \ =\ \ 
N\cap \mathrm{aff}(P_1\cup P_2).
\]
\label{JOINS}
\vskip-4mm
\end{theorem}

\noindent \textit{Proof.} By \textsf{(e)}  it is obvious  that
$\mathcal{T}_1*\mathcal{T}_2$ is a lattice  triangulation of $P$\/ w.r.t.~the
lattice $N$.
\par\noindent
(i) Let $\mathbf{s}$ be a $(d_1+d_2+1)$-simplex of
$\mathcal{T}_1*\mathcal{T}_2$.   
By definition it can be written as
$$
\mathbf{s}=\mathrm{conv}
(\{n^{(1)}_0,\dots,n^{(1)}_{d_1},n^{(2)}_0,\dots,n^{(2)}_{d_2}\}),
$$
where
$\text{\rm conv}(\{n^{(i)}_0,\dots,n^{(i)}_{d_i}\})\in   \mathcal{T}_i(d_i)$  
is a basic simplex w.r.t.~$N\cap \text{\rm aff}(P_i)$, for       $i=1,2$. In 
particular, we have 
$\mathrm{aff}_\Bbb{Z}(\{n^{(1)}_0,\dots,n^{(1)}_{d_1},n^{(2)}_0,\dots,n^{(2)}_{d_2
}\})=
\mathrm{aff}_\Bbb{Z}(\mathrm{vert}(\mathcal{T}_1)\cup 
\mathrm{vert}(\mathcal{T}_2))$. 
Thus, $\mathbf{s}$ is basic w.r.t.~$N\cap\mathrm{aff}(P_1\cup P_2)$
if and only if 
\[
\mathrm{aff}_{\Bbb{Z}}\big(\,
\mathrm{vert}(\mathcal{T}_1)\cup 
\mathrm{vert}(\mathcal{T}_2)\,\big) \ \ =\ \ 
N\cap \mathrm{aff}(P_1\cup P_2).
\]
\par\noindent (ii)  If $\psi _1\in  $ SUCSF$_{\Bbb{R}}\left( \mathcal{T}_1\right)  $,
$\psi _2\in $ SUCSF$_{\Bbb{R}}\left( \mathcal{T}_2\right)  $, then the function $\Psi
$ defined~by
\[
\Psi \left( t\ \mathbf{x}+\left( 1-t\right) \mathbf{\ y}\right) :=t\ \psi
_1\left( \mathbf{x}\right) +\left( 1-t\right) \ \psi _2\left( \mathbf{y}\right) 
,\text{ for all \ }\mathbf{x}\in P_1,\,\mathbf{y}\in P_2,\ \text{and\ }\emph{\ 
}t\in \left[ 0,1\right] ,
\]
belongs to SUCSF$_{\Bbb{R}}\left( \mathcal{T}_1*\mathcal{T}_2\right) $.
\par\noindent (iii) If $\mathcal{T}_1$ and $\mathcal{T}_2$ are both balanced
and have dimensions $d_1$ and $d_2$, with colouring functions
\[
\varphi _1:\text{vert}\left( \mathcal{T}_1\right) \rightarrow \left\{ 0,1,2,\ldots
,d_1\right\} \ \ \ \text{ and \ \ \ }\varphi _2:\text{vert}\left( 
\mathcal{T}_2\right)
\rightarrow \left\{ 0,1,2,\ldots ,d_2\right\},
\]
then $\varphi _1*\varphi _2:$ vert$\left( \mathcal{T}_1*\mathcal{T}_2\right) =$
vert$\left( \mathcal{T}_1\right) \cup $ vert$\left( \mathcal{T}_2\right) 
\rightarrow \left\{
0,1,\ldots ,d_1+d_2+1\right\} $ defined~by
\[
\left( \varphi _1*\varphi _2\right) \left( \mathbf{x}\right) :=\left\{
\begin{array}{ll}
\varphi _1\left( \mathbf{x}\right)  , & \text{if \ \ }\mathbf{x\in }\text{
vert}\left( \mathcal{T}_1\right) \\
\varphi _2\left( \mathbf{x}\right) +d_1+1  , & \text{if \ \ }\mathbf{x\in }\text{ 
vert}\left( \mathcal{T}_2\right)
\end{array}
\right.
\]
is   a  colouring   function  for   $\mathcal{T}_1*\mathcal{T}_2$.  $_{\Box}$
\eject

\noindent \textsf{(f) }Let now $\mathcal{T}$ be a triangulation of $d$-polytope 
$P\subset \Bbb{R}^d$. If $\lambda \,P=\left\{ \lambda \,\mathbf{x\
}\left| \mathbf{x}\in P\right. \right\} $ is the $\lambda $-times \textit{dilated 
}$P$, i.e., the image of $P$ under the dilation map 
$d_{\lambda}:\Bbb{R}^d\ni
\mathbf{x}\mapsto \lambda \,\mathbf{x}\in \Bbb{R}^d$, where $\lambda \in
\Bbb{Z}$, $\lambda \geq 1$, then $\lambda \,\mathcal{T}$ will denote the
corresponding triangulation of $\lambda \,P$ under the same dilation map,
i.e., $\lambda \,\mathcal{T}=\left\{ \lambda \,\mathbf{s\ }\left| \ \mathbf{s}\in 
\mathcal{T}\right. \right\} $.

\begin{lemma}[Preservation of the ``c.b.-property'' for dilations]
\label{dilations}\ \ \newline
Let $\mathcal{T}$ denote a triangulation of a $d$-polytope $P\subset \Bbb{R}^d$, 
and $\lambda\geq 1$ an integer.\smallskip \newline
\emph{(i) }If\/ $\mathcal{T}$ is coherent, then $\lambda \,\mathcal{T}$ is
coherent.\smallskip \newline
\emph{(ii) }If\/ $\mathcal{T}$ is balanced, then $\lambda \,\mathcal{T}$ is
balanced.\smallskip \newline
\emph{(iii) }If $P$ is a lattice polytope w.r.t.\ a lattice $N\subset
\Bbb{R}^d$, and $\mathcal{T}$ a lattice triangulation, then $\lambda 
\,\mathcal{T}$ is a lattice triangulation too w.r.t.\ the same lattice.
\end{lemma}

\noindent \textit{Proof. }For (i), if $\psi \in $ SUCSF$_{\Bbb{R}}( \mathcal{T})$, 
then $\psi\circ (d_\lambda)^{-1} \in $ SUCSF$_{\Bbb{R}}( \lambda \,\mathcal{T})$. 
For (ii), if $\varphi $ is a colouring function for $\mathcal{T}
$, then $\varphi\circ(d_\lambda)^{-1}$ serves as a colouring function for $\lambda 
\,\mathcal{T}$. (iii) is obvious because $\lambda $ is an integer. 
$_{\Box }$

\begin{remark}
\emph{If} $\mathcal{T}$ \emph{is a basic triangulation of a lattice} 
$d$\emph{-polytope} $P\subset \Bbb{R}^d$\emph{, then for} $\lambda
\geq 2$\emph{, the dilation }$\lambda \,\mathcal{T}$ \emph{is non-basic (see
Figure 2). Nevertheless, as} \emph{we shall see below in 
Theorem~\ref{Watti}, any
dilation of a b.c.b.\ triangulation admits a b.c.b.\ triangulation.}
\end{remark}

\[
\unitlength 0.6truecm
\begin{picture}(10,5)
\multiput(0,0)(1,0){3}{\circle*{0.15}}
\multiput(1,1)(1,0){2}{\circle*{0.15}}
\multiput(2,2)(1,0){1}{\circle*{0.15}}
\thicklines
\put(0,0){\line(1,0){2}}
\put(1,1){\line(1,0){1}}
\put(2,0){\line(0,1){2}}
\put(1,0){\line(0,1){1}}
\put(0,0){\line(1,1){2}}
\put(1,0){\line(1,1){1}}
\put(1,-0.8){$\mathbf{\mathcal{T}}$}
\put(8,0){%
\multiput(0,0)(1,0){7}{\circle*{0.15}}
\multiput(1,1)(1,0){6}{\circle*{0.15}}
\multiput(2,2)(1,0){5}{\circle*{0.15}}
\multiput(3,3)(1,0){4}{\circle*{0.15}}
\multiput(4,4)(1,0){3}{\circle*{0.15}}
\multiput(5,5)(1,0){2}{\circle*{0.15}}
\multiput(6,6)(1,0){1}{\circle*{0.15}}
\thicklines
\put(0,0){\line(1,0){6}}
\put(3,3){\line(1,0){3}}
\put(6,0){\line(0,1){6}}
\put(3,0){\line(0,1){3}}
\put(0,0){\line(1,1){6}}
\put(3,0){\line(1,1){3}}
\put(3,-0.8){$\mathbf{3\mathcal{T}}$}
}
\end{picture}
\]
\begin{center} {\bf Figure 2} \end{center}

\section{Torus-equivariant, crepant, projective resolutions\label{TOR-CR-R}}

\noindent Reverting to (singular) Gorenstein abelian quotient spaces 
$$
X\left( N_G,\Delta _G\right) =U_{\sigma _0}=\Bbb{C}^d/G,\quad G\subset  \hbox{\rm 
SL}\left( d,\Bbb{C}\right),\quad d\geq 2,
$$
we explain how the desired
desingularizations can be constructed by means of triangulations.

\noindent Let $f=$ id$_{*}:X( N_G,\widehat{\Delta _G})
\rightarrow X\left( N_G,\Delta _G\right) $ be an arbitrary $T_{N_G}$-equivariant 
desingularization of $X\left( N_G,\Delta _G\right) $ 
(as in Thm.~\ref{TORRES}). There are one-to-one correspondences:

\begin{equation}
\begin{array}{ccc}
\varrho & \in & \left( \widehat{\Delta _G}\left( 1\right) \smallsetminus
\left\{ \Bbb{R}_{\geq 0}\,e_1,\ldots ,\Bbb{R}_{\geq 0}\,e_d\right\} \right)
\\
\updownarrow &  & \updownarrow \\
n\left( \varrho \right) & \in & \left( \bigcup \left\{ \text{Sk}^1\left(
\sigma \right) \,\left| \,\sigma \in \widehat{\Delta _G}\right. \right\}
\smallsetminus \left\{ e_1,\ldots ,e_d\right\} \right) \\
\updownarrow &  & \updownarrow \\
D_{n\left( \varrho \right) }:=V\left( \varrho \right) =V\left( \Bbb{R}_{\geq
0}\,n\left( \varrho \right) \right) & \in & \left\{
\begin{array}{c}
\text{exceptional prime divisors} \\
\text{with respect to }f
\end{array}
\right\}
\end{array}
\label{exc-div}
\end{equation}
Furthermore, $D_{e_i}:=V( \Bbb{R}_{\geq 0}\,e_i) $, in $X(
N_G,\widehat{\Delta _G}) $, corresponds to the strict transform of
\[
\left( \left\{ \mathbf{z}=\left( z_1,\ldots ,z_d\right) \in \Bbb{C}^d\,\left| 
\,z_i=0\right. \right\} \right) /G
\]
with respect to $f$, for $1\leq i\leq d$.

\begin{proposition}
A\/ $T_{N_G}$-equivariant partial desingularization $f:X\left( N_G,\widehat{\Delta 
_G}\right) \rightarrow U_{\sigma _0}$
\emph{(}defined by any fan $\widehat{\Delta _G}$ refining $\Delta _G$ w.r.t.\ 
$N_G\emph{)}$ is crepant if
and only if
\begin{equation}
\bigcup \left\{ \text{\emph{Sk}}^1\left( \sigma \right) \,\left| \,\sigma
\in \widehat{\Delta _G}\right. \right\} \subset \left\{ \mathbf{x=}\left(
x_1,\ldots ,x_d\right) ^{\intercal }\in \left( N_G\right) _{\Bbb{R}}\ \left|
\ \sum_{i=1}^d\ x_i=1\right. \right\} .  \label{crep-prop}
\end{equation}
\end{proposition}

\noindent \textit{Proof. }Let $m_1,\ldots ,m_d$ be $\Bbb{R}$-linearly
independent vectors of $M_G$. The dualizing sheaf $\omega _{T_{N_G}}$ of $T_{N_G}$ 
is generated by the rational differential form
\[
\phi _{T_{N_G}}:=\frac{\pm\text{det}(M_G)}
{\text{det}\Big(\bigoplus\limits_{1\le i\le d}\Bbb{Z} m_i\Big)}
\cdot \frac{du_1}{u_1}\wedge \cdots \wedge \frac{du_d}{u_d}\,\,,
\]
where $u_1=\mathbf{e}\left( m_1\right) ,\ldots ,u_d=\mathbf{e}\left(
m_d\right) $ ( $\phi _{T_{N_G}}$ is independent of the specific choice of 
$m_1,\ldots ,m_d$).\newline
Moreover, for the dualizing sheaf $\omega _{X\left( N_G,\Delta _G\right) }$
of $X\left( N_G,\Delta _G\right) $ we have
\[
\omega _{X\left( N_G,\Delta _G\right) }=\mathcal{O}_{X\left( N_G,\Delta
_G\right) }\left( K_{X\left( N_G,\Delta _G\right) }\right) =\Bbb{C}\left[
M_G\cap \text{ int}\left( \sigma _0^{\vee }\right) \right] \cdot \phi
_{T_{N_G}}
\]
(cf.\ \cite[Lemma 3.3, p.~293]{Reid1}). On the other hand, conditions (i)-(v)
of Prop.~\ref{Gor-prop} are equivalent to the triviality of $\omega
_{X\left( N_G,\Delta _G\right) }$, as one can easily verify by using
Ishida's criteria \cite[p.~126]{Oda}. This means, in particular, that
the semigroup ideal $M_G\cap $ int$\left( \sigma _0^{\vee }\right) $ of the
semigroup ring $M_G\cap\sigma _0^{\vee }$ is principal. In fact, it is
generated by the element $( 1,1,\ldots ,1,1) $ and $\Bbb{C}\left[
M_G\cap \text{ int}\left( \sigma _0^{\vee }\right) \right] $ by $\mathbf{e}\left( 
\left( 1,1,\ldots ,1,1\right) \right)$,
where $\text{ord}_{V(\Bbb{R}_{\ge0}e_i)} \big(\mathbf{e}((1,\ldots,1))\cdot
\phi_{T_{N_G}}\big)=1$ for $1\le i\le d$.
\smallskip

$f$ is crepant if and only if the difference
\[
K_{X\left( N_G,\widehat{\Delta _G}\right) }-f^{*}\left( K_{X\left(
N_G,\Delta _G\right) }\right)
\]
between the canonical divisor of $X\left( N_G,\widehat{\Delta _G}\right) $
and the pull-back of the canonical (trivial) divisor of $X\left( N_G,\Delta
_G\right) $ vanishes. In this case, $K_{X\left( N_G,\widehat{\Delta _G}\right) }$ 
is trivial too. Let
\[
X( N_G,\widehat{\Delta _G}) =\Big( \coprod_{\sigma \in
\widehat{\Delta _G}}\,\widehat{U}_\sigma \Big) \Big/\sim ,
\]
$\sigma '\in \widehat{\Delta _G}\left( d\right) $ an arbitrary
maximal cone of $\widehat{\Delta _G}$, $\varrho $ a ray from $\left(
\widehat{\Delta _G}\left( 1\right) \smallsetminus \left\{ \Bbb{R}_{\geq
0}\,e_1,\ldots ,\Bbb{R}_{\geq 0}\,e_d\right\} \right) $ with $\varrho \prec
\sigma '$, and $D_{n\left( \varrho \right) }$ the corresponding
exceptional prime divisor (as in (\ref{exc-div})). Since
\[
H^0\left( \widehat{U}_{\sigma '};\omega _{X\left( N_G,\widehat{\Delta _G}\right) 
}\right) =\Bbb{C}\left[ M_G\cap \text{ int}\left( \left(
\sigma '\right) _{}^{\vee }\right) \right] \cdot \phi
_{T_{N_G}}\,\,,
\]
by the same argument as above, if $\omega _{\widehat{U}_{\sigma '}}$
is trivial, then it is generated by $\mathbf{e}\left( \left( 1,1,\ldots
,1,1\right) \right) \cdot \phi _{T_{N_G}}$. 
The discrepancy of $D_{n(\varrho)}$ with respect to 
$f\big| _{\widehat{U}_{\sigma '}}$ is equal to 
$\text{ord}_{D_{n(\varrho)}}
\big( \mathbf{e}((1,\ldots ,1))\cdot\phi_{T_{N_G}}\big)-1$,
where this vanishing order
of $\left( \mathbf{e}\left( \left(
1,\ldots ,1\right) \right) \cdot \phi _{T_{N_G}}\right)$
along ${D_{n(\varrho ) }}$ equals
\[
\text{ord}_{D_{n(\varrho)}}
\big( \mathbf{e}((1,\ldots ,1))\cdot\phi_{T_{N_G}}\big)
\ =\ 
\left\langle (1,\ldots,1) , n(\varrho)\right \rangle,
\]
by \cite[lemma of p.~61]{Fulton} and
the definition of $\phi _{T_{N_G}}$. Hence $f$ is crepant if and only if the
condition (\ref{crep-prop}) is satisfied. $_{\Box }$

\begin{corollary}
All $T_{N_G}$-equivariant partial crepant desingularizations of $U_{\sigma
_0}=\Bbb{C}^d/G$ {\rm (}with overlying spaces having at most Gorenstein abelian 
quotient singularities{\rm)} are of the form
\begin{equation}
f_{\mathcal{T}}:X\left( N_G,\widehat{\Delta _G}\left( \mathcal{T}\right)
\right) \rightarrow U_{\sigma _0}  \label{CRDES}
\end{equation}
for some fan $\widehat{\Delta _G}\left( \mathcal{T}\right) =\left\{ \sigma 
_{\mathbf{s}},\mathbf{s}\in \mathcal{T}\right\} $, where $\sigma 
_{\mathbf{s}}:=\left\{ \gamma \mathbf{x\,\in }\left( N_G\right) _{\Bbb{R}}\,\left| 
\,\mathbf{x}\in \mathbf{s},\gamma \in \Bbb{R}_{\geq 0}\,\right. \right\} $ and
$\mathcal{T}$ denotes a lattice triangulation of the junior simplex $\frak{s}_G$.
\end{corollary}

\noindent 
Since a simplex $\mathbf{s}\in \mathcal{T}$ is basic if and only
if mult$\left( \sigma _{\mathbf{s}}, N_G\right) =1$, i.e., if and only if the
toric variety~$U_{\sigma _{\mathbf{s}}}$ is smooth, we have:

\begin{lemma}
$f_{\mathcal{T}}$ $:X\left( N_G,\widehat{\Delta _G}\left( \mathcal{T}\right)
\right) \rightarrow U_{\sigma _0}$ is a crepant \emph{(}full\emph{)}
desingularization if and only if $\mathcal{T}$ is basic.
\end{lemma}

\begin{example}
\emph{Consider the }$3$\emph{-dim.~Gorenstein cyclic quotient
singularity of} \textit{type} $\frac 17(3,3,\allowbreak 1) $ \emph{(in the
sense of \cite[\S\ 4.2]{Reid3}). The junior simplex} $\frak{s}_G=$
\emph{conv}$\left( \left\{ e_1,e_2,e_3\right\} \right) $ \emph{contains three
additional inner lattice points, namely} $n_1=\frac 17\left( 3,3,1\right)
^{\intercal }$\emph{,} $n_2=\frac 17\left( 2,2,3\right) ^{\intercal }$\emph{,} 
$n_3=\frac 17\left( 1,1,5\right) ^{\intercal }$\emph{. The unique} \emph{basic 
triangulation} $\mathcal{T}$ \emph{of} $\frak{s}_G$ \emph{is drawn in
Figure 3.} \emph{The three exceptional prime divisors on} $X\left( 
N_G,\widehat{\Delta _G}\left( \mathcal{T}\right) \right) $ 
\emph{are isomorphic to}
\[
D_{n_1}\cong \Bbb{P}_{\Bbb{C}}^2\,,\ \ \ D_{n_2}\cong \Bbb{F}_3,\text{ \ \
\emph{and} \ \ }D_{n_3}\cong \Bbb{F}_5,
\]
\emph{respectively, where }$\Bbb{F}_j:=\Bbb{P}\left( 
\mathcal{O}_{\Bbb{P}_{\Bbb{C}}^1}\oplus \mathcal{O}_{\Bbb{P}_{\Bbb{C}}^1}\left( 
j\right) \right)
,j\geq 0$, \emph{denote the Hirzebruch-surfaces.}
\end{example}

\[
\unitlength1truecm
\begin{picture}(0,0)
\put(5,2.1){$f_{\mathcal{T}}$}
\put(-.4,0){$e_1$}
\put(3.9,0){$e_2$}
\put(2,3.3){$e_3$}
\put(5.8,0){$e_1$}
\put(10,0){$e_2$}
\put(8.2,3.3){$e_3$}
\put(2,0.4){$n_1$}
\put(2,1.2){$n_2$}
\put(2,2){$n_3$}
\end{picture}
\def\epsfsize#1#2{.7#1}
\epsffile{hirzebruch.eps}
\]
\begin{center} {\bf Figure 3} \end{center}

\noindent Coming back to (\ref{CRDES}), we see that for every $\psi\in
\mathrm{SUCSF}_{\Bbb{Q}}
\left( N_G,\widehat{\Delta _G} ( \mathcal{T})\right)$, the restriction $\psi
| _{\mathcal{T}}$ belongs to $\mathrm{SUCSF}_{\Bbb{R}}( \mathcal{T})$; and conversely, 
starting from any $\psi\in \mathrm{SUCSF}_{\Bbb{R}}( \mathcal{T})$, one can construct by 
``pulling vertices'' a strictly upper convex $\mathcal{T}$-support function $\psi^\prime$
with $\psi^\prime(\mathrm{vert}(\mathcal{T}))$ belonging to $\Bbb{Q}$ (or even to $\Bbb{Z}$ after appropriate ``scaling''), and take
an extension $\widetilde{\psi}$ of it on $| \widehat{\Delta _G} (
\mathcal{T})|$, for which $\widetilde{\psi }\left( \gamma
\mathbf{x}\right) =\psi^\prime \left( \mathbf{x}\right) $, for
all $\mathbf{x}\in \mathbf{s}$, $\mathbf{s}\in \mathcal{T}\left( d\right) $,
and $\gamma \in \Bbb{R}_{\geq 0}$). Therefore, by \ref{SUCF}, \ref{quasiproj}, we 
obtain:

\begin{proposition}
$f_{\mathcal{T}}$ is a projective crepant morphism if and only if $\mathcal{T}$ is 
coherent.
\end{proposition}

\noindent Hence, the existence of a projective, crepant (full)
desingularization of $U_{\sigma _0}=\Bbb{C}^d/G$ is equivalent to the
existence of a basic, coherent triangulation of the junior simplex $\frak{s}_G$. 
Our use of balanced complexes is a requirement 
which is not needed from the toric point of view; however,
it is our main technical tool for
``gluing'' triangulations together.
In particular, it is not clear whether every dilation of
a basic triangulation has a basic refinement; but this will 
be proved for balanced basic triangulations, and for
b.c.b.\ triangulations, below. 

\section{From special data to Watanabe simplices}

\noindent In 1980 Watanabe \cite[1.7, 2.1 and 2.8]{Wata} classified
\textit{all} abelian quotient c.i.-singularities 
in terms of ``special data'' which encode 
complete descriptions of the corresponding polynomial rings and acting
groups. Although his method for proving this theorem is purely algebraic, he
also indicated how could one represent these ``special data'' by certain
graphs (which are actually forests, and will be henceforth called \textit{Watanabe 
forests}). As we shall see in the present section, this
graph-theoretic approach to the Classification Theorem \ref{WAT-TH} 
yields a very useful lattice-geometric interpretation. More precisely, we work
out the following one-to-one correspondences:
\[
\fbox{$
\begin{array}{ccc}
&  &  \\
&
\begin{array}{ccc}
\left\{
\begin{array}{c}
\text{special non-trivial} \\
\text{data \ }\Bbb{D=}\left( \frak{D},\mathbf{w}\right) \\
\text{(up to isomorphisms }\simeq \text{) }
\end{array}
\right\} & \stackrel{}{\leftrightarrow } & \left\{
\begin{array}{c}
\text{abelian quotient c.i.-singularities} \\
\left( \Bbb{C}^d/G,\left[ \mathbf{0}\right] \right) \cong \left( 
\Bbb{C}^d/G_{\Bbb{D}},\left[ \mathbf{0}\right] \right) \\
\text{(up to analytic isomorphisms }\cong \text{)}
\end{array}
\right\} \\
\updownarrow &  & \updownarrow \  \\
\left\{
\begin{array}{c}
\begin{array}{c}
\text{non-trivial} \\
\text{Watanabe forests }\mathbf{W}_{\Bbb{D}}
\end{array}
\\
\text{(up to graph-theoretic} \\
\text{isomorphisms }\cong _{\,\text{wg}}\text{)}
\end{array}
\right\} & \stackrel{}{\leftrightarrow } & \left\{
\begin{array}{c}
\text{Gorenstein, simplicial, singular, toric\ } \\
\text{ varieties\ }X\left( N_G,\Delta _G\right) =U_{\sigma _0}\text{ for which } 
\\
\frak{s}_{G\text{ }}\ \text{is a Watanabe simplex w.r.t.\ $N_G$} \\
\text{(up to affine integral transformations)}
\end{array}
\right\}
\end{array}
&  \\
&  &
\end{array}
$}
\]
Let us first formulate Watanabe's result.

\begin{definition}
\label{SP-DATUM}\emph{Let }$d\geq 2$ \emph{be an integer}\emph{. A}
\textit{special datum} $\Bbb{D=}\left( \frak{D},\mathbf{w}\right) $ \emph{(w.r.t.\ 
}$d$\emph{) is a pair consisting of a set of non-empty
subsets of }$\left\{ 1,2,\ldots ,d\right\} $ \emph{(i.e.} $\frak{D}\subseteq
\mathbf{2}^{\{ 1,2,\ldots ,d\} }\smallsetminus\{\varnothing\}$\emph{),} 
\emph{together with a
``weight-function'' }$\mathbf{w}:\frak{D}\rightarrow \Bbb{N}$\emph{, such
that:\smallskip }\newline
\emph{(i) For each } $i\in \left\{ 1,2,\ldots ,d\right\} $
\emph{we have} $\left\{ i\right\} \in \frak{D}$\emph{.\smallskip }\newline
\emph{(ii) For every pair of index-subsets }$J$\emph{,} $J'\in
\frak{D}$\emph{, either }$J\subseteq J'$\emph{,} \emph{or }$J^{\prime
}\subseteq J$\emph{,} \emph{or} $J\cap J'=\varnothing $\emph{.\smallskip }
\newline
\emph{(iii) If }$J$ \emph{is a maximal element with respect to the inclusion
relation ``}$\subseteq$\emph{'', then }$\mathbf{w}\left( J\right) =1$.
\smallskip \newline
\emph{(iv) If }$J$\emph{,} $J'\in \frak{D}$ \emph{and} $J\subsetneqq J'$\emph{, 
then} $\mathbf{w}\left( J\right) >\mathbf{w}\left( J'\right) $ \emph{and} 
$\mathbf{w}\left( J'\right)
\left| \mathbf{w}\left( J\right) \right. $.
\smallskip \newline
\emph{(v) For} $J_1,J_2,J\in \frak{D}$\emph{, } \emph{with} $J_1\sqsubset 
J$\emph{,} $J_2\sqsubset J$\emph{, we have }$\mathbf{w}\left( J_1\right) 
=\mathbf{w}\left( J_2\right) $\emph{.\smallskip }\newline
\emph{For this the
binary cover relation ``}$\sqsubset $\emph{'' on the sets of }$\frak{D}$\emph{\ is
defined by:}
\[
J\sqsubset J'\stackunder{\text{\emph{def}}}{\Longleftrightarrow }\left(
J\subsetneqq J' \mbox{\rm\ \ and \ } 
\nexists J^{\prime \prime },\ J^{\prime \prime }\in \frak{D}:\
J\subsetneqq J^{\prime \prime }\subsetneqq J'\right).
\]
$\Bbb{D=}\left( \frak{D},\mathbf{w}\right) $ \emph{is} \textit{non-trivial 
}\emph{if there is at least one subset} $J\in \frak{D}$ 
\emph{ with} \emph{cardinality} $\#\left( J\right) \geq 2$\emph{.}
\end{definition}

\begin{definition}
\label{DEFGD}\emph{Let }$d\geq 2$ \emph{be an integer, and let } 
$\Bbb{D=}\left( \frak{D},\mathbf{w}\right) $\emph{\ be a special datum w.r.t.\ 
}$d$\emph{. We define the polynomial ring }$R_{\Bbb{D}}$\emph{\ and the group} 
$G_{\Bbb{D}}$ \emph{by}
\[
R_{\Bbb{D}}:=\Bbb{C}\left[ \frak{x}_J\ \left| \ J\in \frak{D}\right. \right]
\text{ \ \ \ \emph{with} \ \ \ }\frak{x}_J:=\Big( \prod_{i\in J}\ \frak{x}_i\Big) 
^{\mathbf{w}( J) }, \ \ \ \mbox{\rm and}
\]
\[
G_{\Bbb{D}}:=\left\langle    \left\{    \text{\emph{diag}}\Big(   1,\ldots
,1,\stackunder{i\text{\emph{-th                   pos.}}}{\underbrace{\zeta
_{\mathbf{w}}}},1,\ldots ,1,\stackunder{j\text{\emph{-th pos.}}}{\underbrace{\zeta 
_{\mathbf{w}}^{-1}}},1,\ldots,1\Big) \ \left|
\begin{array}{c}
i\in J_1,\ j\in J_2 \\
\text{\emph{for all \ }}J_1,J_2,J\in \frak{D} \\
\text{\emph{with \ }}J_1\sqsubset J\emph{,\ }J_2\sqsubset J \\
\text{\emph{and \ }}\mathbf{w}=\mathbf{w}\left( J_1\right) =\mathbf{w}\left(
J_2\right)
\end{array}
\right. \right\} \right\rangle
\]
\emph{Here $\zeta _{\mathbf{w}}$ denotes a
primitive $\mathbf{w}$th root of unity,
and the diagonal matrices generating }$G_{\Bbb{D}}$
\emph{are }$\left( d\times d\right) $\emph{-matrices in} \emph{SL}$\left(
d,\Bbb{C}\right) $.
\end{definition}

\begin{theorem}[Watanabe's Classification Theorem]
\label{WAT-TH}Let $d\ge 2$ be an integer and $G$ be a finite, abelian
subgroup of \emph{SL}$\left( d,\Bbb{C}\right) $. The
following conditions are equivalent\emph{:\smallskip }\newline
\emph{(i) }The quotient space $\Bbb{C}^d/G$ is minimally embeddable as a
complete intersection of hypersurfaces in an affine \emph{(}complex\emph{)}
space.\smallskip \newline
\emph{(ii) }There is a special datum $\Bbb{D=}\left( \frak{D},\mathbf{w}\right) $ 
\emph{(}w.r.t.\ $d$\emph{)}, such that
\[
\Bbb{C}^d/G\cong \text{\emph{Max-Spec}}\left( R_{\Bbb{D}}\right)
\]
and $G$ is conjugate to $G_{\Bbb{D}}$ \emph{(}within \emph{SL}$\left( 
d,\Bbb{C}\right) $\emph{)}.\smallskip 

In other words,
\[
\left( \Bbb{C}^d/G,\left[ \mathbf{0}\right] \right) \cong \left( 
\Bbb{C}^d/G_{\Bbb{D}},\left[ \mathbf{0}\right] \right) ,
\]
i.e., up to an analytic isomorphism, all \emph{(}germs of\emph{)} abelian
quotient \emph{c.i.-}singularities are pa\-ram\-e\-ter\-ized by the set of
non-trivial special data $\Bbb{D=}\left( \frak{D},\mathbf{w}\right) $ 
\emph{(}w.r.t.\ $d$\emph{).\ }
\end{theorem}

\begin{remark}
\label{MAXE}\emph{Let} $\Bbb{D=}\left( \frak{D},\mathbf{w}\right) $ \emph{be
a special datum (w.r.t.} $d$\emph{) and let }$J_{\blacklozenge }^{\left[
1\right] }$\emph{,}$\ldots $\emph{, }$J_{\blacklozenge }^{\left[ \kappa
\right] }$ \emph{denote the maximal elements of} $\frak{D}$ \emph{(w.r.t.\ 
``}$\subseteq$\emph{''). By the properties (i) and (ii) in \ref{SP-DATUM}, we
have:\smallskip }\newline
\emph{(i) }$\bigcup_{i=1}^\kappa \ J_{\blacklozenge }^{\left[ i\right]
}=\left\{ 1,2,\ldots ,d\right\} $ \emph{(set-theoretically).\smallskip }\newline
\emph{(ii) For any }$J\in \frak{D}$\emph{ that is not a singleton
$(\#(J)\ge2)$, there is a set-theoretic
partition:} $J=\,\stackrel{\bullet }{\bigcup }\ \left\{ J'\in
\frak{D\ }\left| \ J'\sqsubset J\right. \right\} .$
\end{remark}

\begin{remark}
\label{CONV}\emph{The elements of }$J$\emph{'s determine the positions of
the roots of unity within the diagonals of the matrices generating 
}$G_{\Bbb{D}}$\emph{. As there are many choices for the elements of }$J$\emph{'s
leading to the same matrices (up to permutations of their entries), and
since Theorem~\ref{WAT-TH}}$\,$\emph{\ gives the classification }\textit{up to
isomorphism}\emph{, we shall define ``isomorphisms'' between special data in
order to work only with convenient representatives of the
corresponding equivalence-classes. For a special datum }$\Bbb{D=}\left(
\frak{D},\mathbf{w}\right) $ \emph{(w.r.t.\ }$d$\emph{) let}
\[
\frak{D}\left( p\right) :=\left\{ J\in \frak{D\ }\left| \
\#\left( J\right) =p\right. \right\} ,\ \ 1\leq p\leq d,
\]
\emph{denote the subset of }$\frak{D}$ \emph{consisting of all index-sets of
fixed cardinality} $p$.
 
\emph{Two special data}
$\Bbb{D} =(\frak{D} ,\mathbf{w} )$ \emph{and}
$\Bbb{D}'=(\frak{D}',\mathbf{w}')$ \emph{(w.r.t.} $d$\emph{) 
are} \textit{isomorphic} \emph{(and we denote this isomorphism by }$\Bbb{D}\simeq 
\Bbb{D}'$\emph{) if there exists a bijection}
$\Theta :\{1,\ldots,d\}\rightarrow \{1,\ldots,d\}$ 
\emph{for which}
\smallskip 
\newline
\emph{(i) } $J\in \frak{D} \Longleftrightarrow 
   \Theta (J) \in \frak{D}'$\emph{ \ (i.e.,
$\Theta$ induces bijections $\frak{D}(p)\leftrightarrow\frak{D}'(p)$
for} $1\leq p\leq d$\emph{), and}
\smallskip 
\par\noindent
\emph{(ii) }$\mathbf{w}'\left( \Theta \left( J\right) \right) =\mathbf{w}\left( 
J\right) $\emph{, for all} $J\in \frak{D}$.
\smallskip 
\newline
\emph{It  is  easy  to  verify  the  following  equivalence-implications:
}\medskip\par\noindent
\[
\Bbb{D}\simeq \Bbb{D}'\Longleftrightarrow \left(
\begin{array}{c}
G_{\Bbb{D}}\ \text{\emph{and} }G_{\Bbb{D}'}\ \text{\emph{belong to
the same}} \\
\text{\emph{conjugacy class (within SL}}\left( d,\Bbb{C}\right) \text{\emph{)}}
\end{array}
\right) \Longleftrightarrow \left( \Bbb{C}^d/G_{\Bbb{D}},\left[ \mathbf{0}\right] 
\right) \cong \left( \Bbb{C}^d/G_{\Bbb{D}'},\left[ \mathbf{0}\right] \right)
\]
\textit{Convention A.} \emph{From now on, as representatives} $\Bbb{D=}\left( 
\frak{D},\mathbf{w}\right) $ \emph{of the equivalence-classes from}
\[
\left( \left\{ \text{\emph{all special data (w.r.t.} }d\text{\emph{)}}\right\} 
/\simeq \right) ,
\]
\emph{we shall consider (without loss of generality) only }$\frak{D}$\emph{'s} 
\emph{all of whose index subsets} $J$ \emph{have
the form }$J=\left\{ \nu ,\nu +1,\ldots ,\xi -1,\xi \right\} $\emph{,} 
$1\le\nu\le\xi \le d$\emph{, i.e., contiguous segments
of }$\left\{ 1,2,\ldots ,d\right\} $\emph{.\smallskip }\newline
\textit{Convention B.}
\emph{We refer to the minimum and maximum elements of} 
$J=\left\{ \nu ,\nu +1,\ldots ,\xi -1,\xi \right\}$ \emph{with the
notation} $\nu=\nu _J$ \emph{and} $\xi=\xi _J$.
\end{remark}

\begin{definition}
\label{SUBDATUM}\emph{Let }$d\ge2$ \emph{be an integer,}  $\Bbb{D=}\left( 
\frak{D},\mathbf{w}\right) $\emph{\ a special datum w.r.t.\ }$d$\emph{, and }
\[
J_0=\left\{ \nu _0,\nu _0+1,\ldots ,\xi _0-1,\xi _0\right\} \in \frak{D}
\]
\emph{a} \textit{fixed} \emph{index-set. We define the} \textit{subdatum} 
$\Bbb{D}_{J_0}\Bbb{=}\left( \frak{D}_{J_0},\mathbf{w}_{J_0}\right) $ \emph{of}
$\Bbb{D}$ \textit{relative to} $J_0$ \emph{by}
\[
\frak{D}_{J_0}:=\left\{ J\in \frak{D\ }\left| \ J\subseteq J_0\right. \right\}
,\ \ \ \ \text{\emph{and\ \ \ \ }}\mathbf{w}_{J_0}\left( J\right) 
:=\frac{\mathbf{w}\left( J\right) }{\mathbf{w}\left( J_0\right) },\ 
\mbox{\rm for all }\ J\in \frak{D}_{J_0}.
\]
\emph{The subdatum} $\Bbb{D}_{J_0}\Bbb{=}\left( 
\frak{D}_{J_0},\mathbf{w}_{J_0}\right) $ \emph{can be viewed as an ``autonomous'' 
special datum }${\Bbb D}'=( \frak{D}',\mathbf{w}^{\prime
}) $ \emph{w.r.t.} $\xi _0-\nu _0+1$ \emph{via the bijection}
\[
\vartheta _0:\frak{D}_{J_0}\rightarrow \frak{D}'=\left\{ 1,2,\ldots
,\xi _0-\nu _0+1\right\}
\]
\emph{defined by}
\[
\vartheta _0\left( \nu _0\right) =1\emph{,\ }\vartheta _0\left( \nu
_0+1\right) =2\emph{,\ }\ldots \ \emph{,}\vartheta _0\left( \xi _0\right)
=\xi _0-\nu _0+1,
\]
\emph{(where} $\mathbf{w}'\left( J\right) :=\mathbf{w}_{J_0}\left(
\vartheta _0^{-1}\left( J\right) \right) $\emph{,} \emph{for all} $J\in
\frak{D}'$\emph{).\smallskip }\newline
\textit{Convention C. }\emph{From now on we shall use a subdatum of a
special datum }$\Bbb{D}$ \emph{relative to some} $J_0$ \emph{in} \emph{both
``roles'' without referring explicitly to the identification map }$\vartheta
_0$. \emph{(In the first case we shall emphasize the induced embedding} 
$R_{\Bbb{D}_{J_0}}\hookrightarrow R_{\Bbb{D}}$ \emph{or }$G_{\Bbb{D}_{J_0}}\subset 
G_{\Bbb{D}}$\emph{; in the latter, its own right to enjoy}
\emph{properties (i)-(v) of \ref{SP-DATUM}).}
\end{definition}

\noindent 
Next we recall some definitions concerning
\textit{trees}. We shall mostly use standard terminology. 
As usual, a \textit{graph} is determined by the set of its
\textit{vertices }and the set of its \textit{edges}. A vertex $v$ is a
\textit{neighbour }of another vertex $v'$, if $v$ and $v'$
are adjacent in the considered graph. A graph is \textit{connected} if there
is a path between any two vertices of it. A \textit{cycle} in a graph is a
simple path from a vertex to itself. By a \textit{weighted graph} is meant a
graph with weights assigned to its vertices. Two weighted graphs are 
\textit{isomorphic }to each other (denoted by $\cong _{\,\text{wg}}$\thinspace ) 
if
there exists a bijection between the sets of their vertices preserving both
adjacency and weights.\medskip

\noindent A graph having no cycle is \textit{acyclic}. A
\textit{tree} is a connected acyclic graph. (A \textit{trivial 
tree}\emph{\textit{\ }}is a tree consisting of only one vertex). An
arbitrary acyclic graph is called a \textit{forest}. (So all
connected components of a forest are trees). A \textit{leaf }is a
vertex of degree at most $1$, i.e., a vertex being contained in 
at most one edge. A \textit{rooted tree }distinguishes one vertex as its 
\textit{root}\emph{. }If $v$ is a non-root vertex of a rooted tree,
then its \textit{parent} is
the neighbour of $v$ on the path connecting $v$ with its root. The 
\textit{children} of $v$ are
its other neighbours. In this case, the leaves are exactly the vertices
without children.\medskip

\noindent A \textit{rooted plane tree} is
a rooted tree which is embedded in the plane, i.e., a rooted tree which
is endowed with a left-to-right ordering specified for the
children of each vertex. (By a \textit{plane forest} we mean
a forest having only rooted plane trees as connected components).

\begin{definition}
\emph{Let }$d$ \emph{be an integer} $\geq 2$\emph{\ and} $\Bbb{D=}\left(
\frak{D},\mathbf{w}\right) $\emph{\ be a special datum w.r.t.\ }$d$\emph{.
The }\textit{Watanabe forest associated to }$\Bbb{D}$\emph{\ is defined to
be a weighted plane forest }$\mathbf{W}_{\Bbb{D}}$ \emph{whose vertices are
in one-to-one correspondence with the elements of }$\frak{D}$ \emph{and
whose edges are to be drawn as follows: If }$J,J'$ \emph{are two
elements of }$\frak{D}$\emph{, then the corresponding vertices of 
}$\mathbf{W}_{\Bbb{D}}$\emph{, say} $v_J$\emph{, }$v_{J'}$\emph{, should be
joinable by an edge of} $\mathbf{W}_{\Bbb{D}}$ \emph{if and only if either} 
$J\sqsubset J'$ \emph{or} $J'\sqsubset J$\emph{. If} $\mathbf{W}_{\Bbb{D}}$ 
\emph{happens to be a tree, then we shall speak of a}
\textit{Watanabe tree} \emph{(being associated to} $\Bbb{D}$\emph{). For a }$J\in 
\frak{D}$\emph{, we associate the weight} $\mathbf{w}\left( v_J\right) 
:=\mathbf{w}\left( J\right) $ \emph{to} $v_J$\emph{. } \emph{Moreover, for} 
$J$\emph{, }$J'\in \frak{D}$\emph{,} \emph{with }$J\sqsubset J^{\prime
}$\emph{, using} \emph{Watanabe's original} \textit{upside-down convention}\emph{, 
we shall join} $v_J$ \emph{and} $v_{J'}$\emph{\ in such a
way,} \emph{that }$v_{J'}$ \emph{lies}\textit{\ over}\emph{\ }$v_J$. \emph{It is 
therefore useful to regard the edge connecting }$v_J$ \emph{with} $v_{J'}$ 
\emph{as ``directed'' and denote it by} $\overrightarrow{v_J,v_{J'}}$. {\rm (For 
typographical reasons,
in our figures we shall denote the vertices of the forest $\mathbf{W}_{\Bbb{D}}$ 
by \fbox{$v_J, \mathbf{w}(J)$} enclosing also their weights.)}
\end{definition}

\begin{remark}
\label{SPLIT}\emph{(i) By \ref{SP-DATUM} (i) we have }$\stackrel{\bullet }{\bigcup 
}\ \left\{ J\ \in \frak{D\ }\left| \ v_J\text{ \ \emph{leaf}}\right. \right\} 
=\left\{ 1,2,\ldots ,d\right\} $ \emph{(set-theoretically).\smallskip }\newline
\emph{(ii) By definition, }
\[
\Bbb{D}\simeq \Bbb{D}'\Longleftrightarrow \mathbf{W}_{\Bbb{D}}\cong
_{\,\text{\emph{wg}}}\mathbf{W}_{\Bbb{D}'}
\]
\newline
\emph{(iii) Let }$\Bbb{D=}\left( \frak{D},\mathbf{w}\right) $ \emph{be a
special datum and assume that its
Watanabe forest }$\mathbf{W}_{\Bbb{D}}$\emph{\ has connected components}
\[
\mathbf{W}_{\Bbb{D}^{\left[ 1\right] }},\mathbf{W}_{\Bbb{D}^{\left[ 2\right]
}},\ldots ,\mathbf{W}_{\Bbb{D}^{\left[ \kappa -1\right] 
}},\mathbf{W}_{\Bbb{D}^{\left[ \kappa \right] }}
\]
\emph{whose vertices are in one-to-one correspondence with the elements of }
\[
\frak{D}^{\left[ 1\right] },\frak{D}^{\left[ 2\right] },\ldots ,\frak{D}^{\left[ 
\kappa -1\right] },\frak{D}^{\left[ \kappa \right] }
\]
\emph{for a partition}
$\frak{D}=\bigcup_{1\leq i\leq \kappa }^{\bullet }\ \frak{D}^{[ i]}$
\emph{of }$\frak{D}$\emph{. Then each component} \emph{\ 
}$\mathbf{W}_{\Bbb{D}^{\left[ i\right] }}$ \emph{has root} $v_{J_{\blacklozenge 
}^{\left[
i\right] }}$\emph{, with} $J_{\blacklozenge }^{\left[ i\right] }$ \emph{as
in \ref{MAXE}, and }$\Bbb{D}^{\left[ i\right] }:=\left( \frak{D}^{\left[
i\right] },\mathbf{w}\left| _{\frak{D}^{\left[ i\right] }}\right. \right) $
\emph{is to be identified with the subdatum} $\Bbb{D}_{J_{\blacklozenge
}^{\left[ i\right] }}$ \emph{(in the sense of \ref{SUBDATUM}), for} 
$1\leq i\leq \kappa $\emph{. Furthermore, the group }$G_{\Bbb{D}}$
\emph{splits into the direct product}
\[
G_{\Bbb{D}}=G_{\Bbb{D}_1}\times G_{\Bbb{D}_2}\times \cdots \times 
G_{\Bbb{D}_{\kappa -1}}\times G_{\Bbb{D}_\kappa }.
\]
\end{remark}

\begin{remark}
\emph{Let }$\Bbb{D=}\left( \frak{D},\mathbf{w}\right) $ \emph{be a special
datum (w.r.t.\ an integer }$d\geq 2$\emph{). Then we have:\smallskip }\newline
\emph{(i) }$\left( \Bbb{C}^d/G_{\Bbb{D}},\left[ \mathbf{0}\right] \right)
\cong \left( U_{\sigma _0},\text{\emph{orb}}\left( \sigma _0\right) \right) $
\emph{is a singularity if and only if }$\mathbf{W}_{\Bbb{D}}$ \emph{contains
at least one non-trivial tree. \smallskip In this case,}
\[
\text{\emph{splcod}}\left( \text{\emph{orb}}\left( \sigma _0\right)
;U_{\sigma _0}\right) =d-\#\,\left\{ \text{\emph{all trivial trees of \ 
}}\mathbf{W}_{\Bbb{D}}\right\} \ .
\]
\newline
\emph{(ii) }$\left( \Bbb{C}^d/G_{\Bbb{D}},\left[ \mathbf{0}\right] \right) $
\emph{is a }$d$\emph{-dimensional msc-singularity if and only if all
connected components }$\mathbf{W}_{\Bbb{D}^{\left[ i\right] }}$ \emph{of 
}$\mathbf{W}_{\Bbb{D}}$ \emph{are non-trivial trees. 
In this case, for }$d\geq 3$\emph{, all }$G_{\Bbb{D}}$\emph{'s are abelian, 
}\textit{non-cyclic}\emph{\textit{\ }groups.
\smallskip }\newline
\emph{(iii) }$\left( \Bbb{C}^d/G_{\Bbb{D}},\left[ \mathbf{0}\right] \right) $
\emph{is a }$d$-\emph{dimensional} \emph{absolutely unbreakable} 
\emph{msc-singularity if and only if }$\mathbf{W}_{\Bbb{D}}$ \emph{itself is a
non-trivial Watanabe tree. If, in addition, this singularity is a
hypersurface-singularity, then Theorem~\ref{WAT-TH} can be simplified as
follows:}
\end{remark}

\begin{proposition}
Let $G$ be a finite abelian subgroup of \emph{SL}$\left( d,\Bbb{C}\right) $,
$d\geq 2,$ of order $l\geq 2$. Then $\left( X\left( N_G,\Delta _G\right) 
,\text{\emph{orb}}\left( \sigma _0\right) \right) $ is an absolutely
unbreakable hypersurface-msc-singularity if and only if $G$ is conjugate
\emph{(}within \emph{SL}$\left( d,\Bbb{C}\right) $\emph{) }to a group of the
form
\[
G\left( d;k\right) :=\left\langle \left\{ \left. \text{\emph{diag}}\left(
1,1,\ldots,1,\stackunder{i\text{\emph{-th position}}}{\underbrace{\zeta 
_k}},\stackunder{\left( i+1\right) \text{\emph{-position}}}{\underbrace{\zeta
_k^{-1}}},1,\ldots,1,1\right) \ \right| \ 1\leq i\leq d-1\right\} \right\rangle
\;,
\]
\newline
for a $k\geq 2$. In this case, $l=k^{d-1}$, and one has the following
analytic germ isomorphisms
$$
\left( X\left( N_G,\Delta _G\right) ,\text{\emph{orb}}\left( \sigma
_0\right) \right) \cong \left( \Bbb{C}^d/G\left( d;k\right) ,\left[ 
\mathbf{0\,}_{\Bbb{C}^d}\right] \right)
$$
and
\begin{equation}
\left( \Bbb{C}^d/G\left( d;k\right) ,\left[ \mathbf{0\,}_{\Bbb{C}^d}\right] 
\right) \cong \left( \left\{ \left(
z_0,..,z_d\right) \in \Bbb{C}^{d+1}\left| \ z_0^k=\prod_{i=1}^dz_i\right.
\right\} ,\mathbf{0\,}_{\Bbb{C}^{d+1}}\right).  \label{hyp-eq}
\end{equation}
In particular, $G\left( d;k\right) =G_{\Bbb{D}}\cong \left( 
\Bbb{Z}/k\Bbb{Z}\right) ^{d-1}$, where $\Bbb{D}=\left( \frak{D},\mathbf{w}\right) 
$ denotes
the special datum with
\[
\frak{D}=\left\{ \left\{ 1\right\} ,\left\{ 2\right\} ,\ldots ,\left\{
d\right\} ,\left\{ 1,2,\ldots ,d\right\} \right\} ,\ \ \mathbf{w}\left(
\left\{ j\right\} \right) =k,\, 1\leq j\leq d,\ \ \emph{and}\text{
\ }\mathbf{w}\left( \left\{ 1,2,\ldots ,d\right\} \right) =1.
\]
\emph{Its associated Watanabe tree is given
in Figure 4.}
\end{proposition}

{ 
\drawwith{%
 \drawwith{%
  \drawwith{%
   \drawwith{%
    \drawwith{%
     \drawwith{%
      \drawwith{%
      \drawline}
     \drawline}
    \dottedline{2}}
   \dottedline{3}}
  \dottedline{2}}
 \drawline}
\drawline}
\setlength{\GapDepth}{1truecm}
\setlength{\GapWidth}{0.5truecm}
\begin{center}
\begin{bundle}{\fbox{$ v_{\{1,\dots,d\}}, { 1} $ } }
        \chunk{\fbox{$ v_{\{1\}}, { k} $ }}
        \chunk{\fbox{$ v_{\{2\}}, { k} $ }}
        \chunk{$\cdots$}
        \chunk{}
        \chunk{$\cdots$}
        \chunk{\fbox{$ v_{\{d-1\}}, { k} $ }}
        \chunk{\fbox{$ v_{\{d\}}, { k} $ }}
\end{bundle}
\vspace{0.3truecm}
\begin{center} {\bf Figure 4} \end{center}
\end{center}
} 

\begin{definition}
\emph{We shall call the above singularities (\ref{hyp-eq}) }$\left(
d;k\right) $\textit{-hypersurface-sin\-gu\-lar\-ities\emph{,} }\emph{or simply 
}$\left( d;k\right) $\textit{-hypersurfaces.}
\end{definition}

\begin{remark}
\label{REMH}\emph{(i)} \emph{Using the projection map}
\[
\left\{ \left( z_0,z_1,\ldots ,z_d\right) \in \Bbb{C}^{d+1}\ \left| \
z_0^k=\prod_{i=1}^dz_i\right. \right\} \ni \left( z_0,z_1,\ldots ,z_d\right)
\longmapsto \left( z_1,\ldots ,z_d\right) \in \Bbb{C}^d
\]
\emph{one may regard} $\left( \Bbb{C}^d/G\left( d;k\right) ,\left[ 
\mathbf{0\,}_{\Bbb{C}^d}\right] \right) $ \emph{as the total space of a} 
$k$-\emph{sheeted} \emph{covering of} $\left( 
\Bbb{C}^d,\mathbf{0\,}_{\Bbb{C}^d}\right) $\emph{\ having the union of all 
coordinate hyperplanes of }$\Bbb{C}^d$ \emph{as its branching locus.\smallskip 
}\newline
\emph{(ii) Removing the assumption for }$\left( X\left( N_G,\Delta _G\right)
,\text{\emph{orb}}\left( \sigma _0\right) \right) $ \emph{to be absolutely
unbreakable, we obtain a direct product of such
hypersurface-singularities.}
\end{remark}

\noindent 
To present a fairly short proof of our main Theorem~\ref{Main} (by
avoiding to work \textit{simultaneously }with forests \textit{and 
}triangulations), we introduce a new term, under the name
``Watanabe simplex.''

\begin{definition}
\emph{Let} $d\ge 0$ \emph{be an integer and} $N$ \emph{a} $d$\emph{-dimensional 
lattice in} $N_{\Bbb{R}}\cong \Bbb{R}^d$\emph{. The}
\textit{Watanabe simplices} \emph{w.r.t.} $N$ \emph{are the lattice simplices}
$\mathbf{s}$ \emph{(of dimension $\le d$) satisfying}
\[
\mathrm{aff}_{\Bbb{Z}}(\mathbf{s}\cap N)\ \ =\ \ \mathrm{aff}(\mathbf{s})\cap N
\]
\emph{which are defined inductively 
(starting in dimension $0$) in the following manner:\smallskip }\newline
\emph{(i) Every $0$-dimensional lattice simplex}
$\mathbf{s}=\left\{ n\right\} $\emph{,} $n\in N$\emph{, is a Watanabe simplex.}
\smallskip 
\newline
\emph{(ii) A lattice simplex} $\mathbf{s}\subset N_{\Bbb{R}}$
\emph{of dimension} $d'$, $1\le d'\le d$,
\emph{is a Watanabe simplex if and only if} \medskip 
\newline
$\bullet$ \textit{either }$\mathbf{s}=\mathbf{s}_1*\mathbf{s}_2$\emph{,
where} $\mathbf{s}_1$, $\mathbf{s}_2$ \emph{are Watanabe simplices of
dimensions} $d_1$, $d_2\ge 0$ with $d'=d_1+d_2+1$,
\emph{with respect to
sublattices} $N_1\subset  \text{\rm aff}(\mathbf{s}_1)$\emph{,} $N_2\subset
\text{\rm aff}(\mathbf{s}_2)$ \emph{of  }$N$\emph{, such that}
$\text{\rm aff}_{\Bbb{Z}}(\mathbf{s}\cap N)=
 \text{\rm aff}_{\Bbb{Z}}(N_1\cup N_2) $,
\medskip 
\newline
$\bullet$ \textit{or }$\mathbf{s}$ \emph{is a lattice translate
of some dilation} $\lambda \,\mathbf{s}'$\emph{,
where} $\lambda \geq 2$ \emph{is an integer, and} 
$\mathbf{s}'$ \emph{is a \mbox{$d'$-dimensional} Watanabe simplex w.r.t.} $N$.
\smallskip \newline
\emph{(These conditions are mutually exclusive;
with this definition every affine integral transformation
that preserves $N$ also preserves the Watanabe simplices of $N$.)}
\end{definition}

\begin{example}
\emph{The }$2$\emph{-dimensional lattice simplices $\mathbf{(A)}$ and 
$\mathbf{(B)}$ (w.r.t.} $\Bbb{Z}^2$\emph{) of Figure 5 are
Watanabe simplices, whereas $\mathbf{(C)}$ and $\mathbf{(D)}$ are
not} \emph{Watanabe simplices.}
\end{example}

\begin{center}
\unitlength 0.6truecm
\begin{picture}(18,3)
\multiput(0,0)(0,1){3}{\multiput(0,0)(1,0){3}{\circle*{0.15}}}
\thicklines
\put(0,0){\line(1,1){2}}
\put(0,0){\line(1,2){1}}
\put(1,2){\line(1,0){1}}
\put(0.5,-1){$\mathbf{(A)}$}
\put(5,0){
\multiput(0,0)(0,1){3}{\multiput(0,0)(1,0){3}{\circle*{0.15}}}
\thicklines
\put(0,0){\line(1,0){2}}
\put(2,0){\line(0,1){2}}
\put(0,0){\line(1,1){2}}
\put(0.5,-1){$\mathbf{(B)}$}
}
\put(10,0){
\multiput(0,0)(0,1){3}{\multiput(0,0)(1,0){3}{\circle*{0.15}}}
\thicklines
\put(0,0){\line(2,1){2}}
\put(0,0){\line(1,2){1}}
\put(2,1){\line(-1,1){1}}
\put(0.5,-1){$\mathbf{(C)}$}
}
\put(15,0){
\multiput(0,0)(0,1){3}{\multiput(0,0)(1,0){5}{\circle*{0.15}}}
\thicklines
\put(0,0){\line(4,1){4}}
\put(4,1){\line(-2,1){2}}
\put(0,0){\line(1,1){2}}
\put(0.5,-1){$\mathbf{(D)}$}
}
\end{picture}
\vspace{0.8truecm}
\begin{center} {\bf Figure 5} \end{center}
\end{center}

\noindent We are now going to prove the following:

\begin{theorem}[Reduction-Theorem]~\\
\label{RED}Let $d\ge 2$ be an integer, 
$\Bbb{D=}\left( \frak{D},\mathbf{w}\right) $ a non-trivial special
datum \emph{(}w.r.t.\ $d$\emph{), }$\mathbf{W}_{\Bbb{D}}$ the associated
Watanabe forest, and $\left( \Bbb{C}^d/G_{\Bbb{D}},\left[
\mathbf{0}\right] \right) $ the corresponding abelian quotient
c.i.-singularity. If we identify the underlying space $\Bbb{C}^d/G_{\Bbb{D}}$
with the toric variety $X\left( N_{G_{\Bbb{D}}},\Delta _{G_{\Bbb{D}}}\right)
=U_{\sigma _0}$ \emph{(}as in \emph{\S\ \ref{TOR-GL}, \textsf{(i)})}, then
the junior simplex $\frak{s}_{G_{\Bbb{D}}}$ is a 
\emph{(}non-basic\emph{)} Watanabe $(d-1)$-simplex w.r.t.~$N_{G_{\Bbb{D}}}$; and conversely, every 
\emph{(}non-basic\emph{)} Watanabe simplex of dimension $\leq d-1$ 
w.r.t.\ an $N\cong \Bbb{Z}^d$ is \emph{(}up to an affine integral
transformation\emph{)} the junior simplex corresponding to some abelian
quotient c.i.-singularity.
\end{theorem}

\noindent The proof will be done in four steps. \smallskip We begin with the
case in which $\mathbf{W}_{\Bbb{D}}$ is a Watanabe tree.\smallskip \newline
$\bullet $ \textit{First step.} Let $\mathbf{W}_{\Bbb{D}}$ denote\textit{\ }a 
non-trivial Watanabe tree. At first we explain how the labeled weights
naturally lead to ``free parameters''. By \ref{MAXE} (i), $J_{\blacklozenge
}=\left\{ 1,2,\ldots ,d\right\} $ is the maximal element of $\frak{D}$
(w.r.t.\ inclusion ``$\subseteq$''). For any $J\in \frak{D}$, $J\subsetneqq 
J_{\blacklozenge }$, there exists a unique chain of oversets
\[
J_0=J\sqsubset J_1\sqsubset J_2\sqsubset \cdots \sqsubset J_{\rho
-1}\sqsubset J_\rho =J_{\blacklozenge }\ \ \ \text{(with }\rho \geq 1\text{)}
\]
which, in fact, corresponds to the directed path
\[
\left[ \overrightarrow{v_{J_0},v_{J_1}},\ \,\overrightarrow{v_{J_1},v_{J_2}},\ 
\ldots \ ,\ \,\overrightarrow{v_{J_{\rho -2}},v_{J_{\rho -1}}},\ 
\,\overrightarrow{v_{J_{\rho -1}},v_{J_\rho }}\right]
\]
connecting $v_J$ with the root $v_{J_{\blacklozenge }}$ of $\mathbf{W}_{\Bbb{D}}$. 
By \ref{SP-DATUM} (iii)-(iv), the weights of the vertices of 
$\mathbf{W}_{\Bbb{D}}$ can be written in the form
\[
\mathbf{w}\left( v_J\right) =\mathbf{w}\left( J\right) =\left\{
\begin{array}{ll}
1  & \text{if\ \ }J=J_{\blacklozenge } \\
\dprod\limits_{i=0}^{\rho -1}\ k_{J_i,J_{i+1}}  & \text{otherwise}
\end{array}
\right.
\]
for integers $k_{J_0,J_1},k_{J_1,J_2},\ldots ,k_{J_{\rho -1},J_\rho }$
(which are $\geq 2$), assigned to the edges
\[
\overrightarrow{v_{J_0},v_{J_1}},\ \,\overrightarrow{v_{J_1},v_{J_2}},\
\ldots \ ,\ \,\overrightarrow{v_{J_{\rho -2}},v_{J_{\rho -1}}},\ 
\,\overrightarrow{v_{J_{\rho -1}},v_{J_\rho }}\ .
\]

\begin{definition}
\emph{By the above procedure we assign an integer} $k_{J,J'}$ $\geq
2$ \emph{to} \textit{every} \emph{edge} $\,\overrightarrow{v_J,v_{J^{\prime
}}}$ \emph{of } $\mathbf{W}_{\Bbb{D}}$ \emph{(}$J\sqsubset J'$\emph{).} 
\emph{These integers} \emph{will be called the} \textit{free parameters} 
\emph{of} \emph{\ }$\Bbb{D}$\emph{\ (or}\textit{\ }\emph{of} 
$\mathbf{W}_{\Bbb{D}}$\emph{). 
For fixed}~$J'$\emph{, the value }$k_{J,J'}$\emph{ is
the }\textit{same }\emph{for} \textit{all} $J$\emph{'s for which
}$J\sqsubset J'$
\emph{\ (see \ref{SP-DATUM} (v)). 
Thus for }$k_{J,J'}$ \emph{we can use the alternative notation}
\[
k_{J,J'}=:k_{\nu _{J'},\xi _{J'}}
\]
\emph{which indicates the ``free parameter source'' (with }$\nu _{J^{\prime
}}$ \emph{and} $\xi _{J'}$ \emph{as defined in \ref{CONV}). For
example, the free parameter }$k$ \emph{of} $\left( d;k\right) 
$\emph{-hypersurfaces equals }$k_{1,d}$.
\end{definition}

\begin{example}
\emph{Up to an analytic isomorphism, all four-dimensional absolutely
unbreakable,} \emph{msc-, quotient c.i.-singularities are }\textit{either}\emph{\ 
}$\left( 4;k\right) $\emph{-hypersurfaces }\textit{or}\emph{\
singularities corresponding to exactly one of the four special data }$\Bbb{D}
$\emph{, whose associated Watanabe trees} $\mathbf{W}_{\Bbb{D}}$ \emph{and
free parameters} $k...$ \emph{are depicted in Figures 6a)-6d).}
\end{example}

\eject

\setlength{\GapDepth}{0.8truecm}
\setlength{\GapWidth}{0.2truecm}
\setlength{\EdgeLabelSep}{0.4truecm}
\begin{center}
\begin{bundle}{\fbox{$ v_{\{1,2,3,4\}}, { 1} $ }}
        \chunk[$k_{1,4}$]{%
        \begin{bundle}{\fbox{$ v_{\{1,2,3\}}, {k_{1,4}} $ } }
        \chunk[$k_{1,3}$]{\fbox{$ v_{\{1\}}, {k_{1,4}\cdot k_{1,3}} $ }  }
        \chunk[$k_{1,3}$]{\fbox{$ v_{\{2\}}, {k_{1,4}\cdot k_{1,3}} $ }  }
        \chunk[$k_{1,3}$]{\fbox{$ v_{\{3\}}, {k_{1,4}\cdot k_{1,3}} $ }  }
        \end{bundle}%
        }
        \chunk[$ k_{1,4}$]{\fbox{$ v_{\{4\}}, {k_{1,4}} $ } }
\end{bundle}
\vspace{0.1truecm}
\begin{center} {\bf Figure 6a)} \end{center}
\end{center}
\medskip

\setlength{\GapDepth}{0.8truecm}
\setlength{\GapWidth}{0.2truecm}
\setlength{\EdgeLabelSep}{0.4truecm}
\begin{center}
\begin{bundle}{\fbox{$ v_{\{1,2,3,4\}}, { 1} $ } }
        \chunk[$k_{1,4}$]{%
        \begin{bundle}{\fbox{$ v_{\{1,2\}}, {k_{1,4}} $ } }
        \chunk[$k_{1,2}$]{\fbox{$ v_{\{1\}}, {k_{1,4}\cdot k_{1,2}} $ }   }
        \chunk[$k_{1,2}$]{\fbox{$ v_{\{2\}}, {k_{1,4}\cdot k_{1,2}} $ }   }
        \end{bundle}%
        }
        \chunk[$ k_{1,4}$]{\fbox{$ v_{\{3\}}, {k_{1,4}} $ } }
        \chunk[$ k_{1,4}$]{\fbox{$ v_{\{4\}}, {k_{1,4}} $ } }
\end{bundle}
\vspace{0.1truecm}
\begin{center} {\bf Figure 6b)} \end{center}
\end{center}
\medskip

\setlength{\GapDepth}{0.8truecm}
\setlength{\GapWidth}{0.2truecm}
\setlength{\EdgeLabelSep}{0.4truecm}
\begin{center}
\begin{bundle}{\fbox{$ v_{\{1,2,3,4\}}, { 1} $ } }
        \chunk[$k_{1,4}$]{%
        \begin{bundle}{\fbox{$ v_{\{1,2\}}, {k_{1,4}} $ }  }
        \chunk[$k_{1,2}$]{\fbox{$ v_{\{1\}},  {k_{1,4}\cdot k_{1,2}}
        $ } }
        \chunk[$k_{1,2}$]{\fbox{$ v_{\{2\}},  {k_{1,4}\cdot k_{1,2}}
        $ } }
        \end{bundle}%
        }
        \chunk[$k_{1,4}$]{%
        \begin{bundle}{\fbox{$ v_{\{3,4\}}, {k_{1,4}} $ }  }
        \chunk[$k_{3,4}$]{\fbox{$ v_{\{3\}},  {k_{1,4}\cdot k_{3,4}}
        $ } }
        \chunk[$k_{3,4}$]{\fbox{$ v_{\{4\}},  {k_{1,4}\cdot k_{3,4}}
        $ } }
        \end{bundle}%
        }
\end{bundle}
\vspace{0.1truecm}
\begin{center} {\bf Figure 6c)} \end{center}
\end{center}
\medskip

\setlength{\GapDepth}{0.8truecm}
\setlength{\GapWidth}{0.2truecm}
\setlength{\EdgeLabelSep}{0.4truecm}
\begin{center}
\begin{bundle}{\fbox{$ v_{\{1,2,3,4\}}, { 1} $ }   }
        \chunk[$k_{1,4}$]{%
        \begin{bundle}{\fbox{$ v_{\{1,2,3\}}, {k_{1,4}} $ }   }
        \chunk[$k_{1,3}$]{%
        \begin{bundle}{\fbox{$ v_{\{1,2\}},  {k_{1,4}\cdot k_{1,3}}
        $ } }
        \chunk[$k_{1,2}$]{\fbox{$     v_{\{1\}},     {k_{1,4}\cdot
        k_{1,3}\cdot k_{1,2}}
        $ } }
        \chunk[$k_{1,2}$]{\fbox{$     v_{\{2\}},     {k_{1,4}\cdot
        k_{1,3}\cdot k_{1,2}}
        $ } }
        \end{bundle}
        }
        \chunk[$k_{1,3}$]{\fbox{$ v_{\{3\}},  {k_{1,4}\cdot k_{1,3}}
        $ } }
        \end{bundle}
        }
        \chunk[$k_{1,4}$]{\fbox{$ v_{\{4\}}, {k_{1,4}} $ }   }
\end{bundle}
\vspace{0.1truecm}
\begin{center} {\bf Figure 6d)} \end{center}
\end{center}
\eject

\noindent $\bullet $ \textit{Second step. }Let $\mathbf{W}_{\Bbb{D}}$ be, as
before, a non-trivial Watanabe tree associated to a special datum $\Bbb{D}=\left( 
\frak{D},\mathbf{w}\right) $ w.r.t.\ a $d\geq 2$. Keeping conventions
A-\thinspace C (of \ref{CONV}, \ref{SUBDATUM}) in mind, we fix an
enumeration, say $J_1,\ldots ,J_p$, of all members of $\left\{ J\in \frak{D}\ 
\left| \ J\sqsubset J_{\blacklozenge }\right. \right\} $, so that
\[
J_i=\left\{ \nu _i,\nu _i+1,\ldots ,\xi _i-1,\xi _i\right\}
\ \ \mbox{\rm for }1\leq
i\leq p,
\]
where
\[
\nu _1=1\leq \xi _1,\ \ \nu _i=\xi _{i-1}+1,
\ \ \mbox{\rm for } 2\leq i\leq p\text{ \ and\ \ }\nu _p\leq \xi _p=d,
\]
and  denote  by  $\mathcal{W}_1,\ldots   ,\mathcal{W}_p$  the  subtrees  of
$\mathbf{W}_{\Bbb{D}}$  which begin  with $v_{J_1},\ldots  ,v_{J_p}$ as  in
Figure 7.

\begin{figure}[htb]
\drawwith{%
 \drawwith{%
  \drawwith{%
   \drawwith{%
    \drawwith{%
     \drawwith{%
      \drawwith{%
      \drawline}
     \drawline}
    \dottedline{2}}
   \dottedline{3}}
  \dottedline{2}}
 \drawline}
\drawline}
\setlength{\GapDepth}{1truecm}
\setlength{\GapWidth}{0.5truecm}
\setlength{\EdgeLabelSep}{0.4truecm}
\begin{center}
\begin{bundle}{\fbox{$ v_{\{1,\dots,d\}}, { 1} $ } }
        \chunk[$k_{1,d}$]{%
        \begin{bundle}{%
        \fbox{$                                  v_{J_1},
        {k_{1,d}} $ }
        }
        \chunk{$\cdots$} $\mathcal{W}_1$
        \chunk{$\cdots$}
        \chunk{$\cdots$}
        \drawwith{\dottedline{3}}
        \end{bundle}
        }
        \chunk[$k_{1,d}$]{
        \begin{bundle}{%
            \fbox{$ v_{J_2}, {k_{1,d}} $ }
        }
        \chunk{$\cdots$} $\mathcal{W}_2$
        \chunk{$\cdots$}
        \chunk{$\cdots$}
        \drawwith{\dottedline{3}}
        \end{bundle}
        }
        \chunk{$\cdots$}
        \chunk{}
        \chunk{$\cdots$}
        \chunk[$k_{1,d}$]{
        \begin{bundle}{
        \fbox{$ v_{J_{p-1}}, {k_{1,d}} $ }
        }
        \chunk{$\cdots$} $\mathcal{W}_{p-1}$
        \chunk{$\cdots$}
        \chunk{$\cdots$}
        \drawwith{\dottedline{3}}
        \end{bundle}
        }
         \chunk[$k_{1,d}$]{
         \begin{bundle}{
         \fbox{$ v_{J_p}, {k_{1,d}} $ }
         }
        \chunk{$\cdots$} $\mathcal{W}_{p}$
        \chunk{$\cdots$}
        \chunk{$\cdots$}
        \drawwith{\dottedline{3}}
        \end{bundle}
        }
\end{bundle}
\vspace{0.3truecm}
\begin{center} {\bf Figure 7} \end{center}
\end{center}
\end{figure}

\noindent To the subtrees $\mathcal{W}_1,\ldots ,\mathcal{W}_p$ we assign
the ``autonomous'' Watanabe trees $\mathbf{W}_{\Bbb{D}_1},\ldots 
,\mathbf{W}_{\Bbb{D}_p}$, induced by the subdata $\Bbb{D}_i=\left( 
\frak{D}_i,\mathbf{w}_i\right) $ of $\Bbb{D}$, where
\[
\frak{D}_i:=\frak{D}_{J_i}\text{,\ \ \ }\mathbf{w}_i:=\mathbf{w}_{J_i},
\ \ \mbox{\rm for } \ \ 1\leq i\leq p,
\]
(in the notation introduced in \ref{SUBDATUM}). $\mathbf{W}_{\Bbb{D}_i}$ is
derived from $\mathbf{W}_{\Bbb{D}}$ after the killing of the parent 
$v_{J_{\blacklozenge }}$ of $v_{J_i}$ (w.r.t.\ $\mathbf{W}_{\Bbb{D}}$), the
christening of $v_{J_i}$ as its root, and the adoption of the appropriate
weights. In other words, the children of the root of $\mathbf{W}_{\Bbb{D}}$
become autonomous roots (and parents).
\medskip 
\newline
$\bullet $ In the arguments of the next four lemmas we shall use induction
on the cardinality number $d$ of $J_{\blacklozenge }$, i.e., on the dimension
of the singularity. (That all assertions are correct for $d=2$ can be
checked easily, and will be therefore omitted). Our ``induction hypothesis''
is that all assertions (formulated below) are true for all
Watanabe trees having roots corresponding to maximal index-sets of
cardinality $<d$, and, in particular, true for $\mathbf{W}_{\Bbb{D}_1},\ldots 
,\mathbf{W}_{\Bbb{D}_p}$. (We do not
exclude the possibility that some of the
$\mathbf{W}_{\Bbb{D}_i}$'s are trivial trees.)

\begin{lemma}
\emph{(i) }If we define
\[
\frak{Y}_{\Bbb{D}}:=
\dbigcup\limits_{J\in\frak{D}}\        
\left\{\text{\emph{diag}}
\Big(   1,\ldots   ,1,
\stackunder{\nu_J\text{\emph{-th pos.}}}
{\underbrace{\zeta    _{\mathbf{w}(    J^{\prime   })}}}
,1,\ldots,1,
\stackunder{\nu_{J'}\text{\emph{-th pos.}}}
{\underbrace{\zeta_{\mathbf{w}(J')}^{-1}}},1,\ldots,1\Big) 
\ \left|
\begin{array}{l}
\text{\emph{for all }}J'\in \frak{D} \\
\text{\emph{with \ }}J'\sqsubset J \\
\text{\emph{and }}\nu _J<\nu _{J'}
\end{array}
\right. \right\} ,\smallskip
\]
then $G_{\Bbb{D}}=\left\langle \frak{Y}_{\Bbb{D}}\right\rangle $, and 
$\frak{Y}_{\Bbb{D}}$ is a minimal generating system for $G_{\Bbb{D}}$.
\smallskip 
\newline
\emph{(ii) }Identifying $\Bbb{C}^d/G_{\Bbb{D}}$ with $X\left( 
N_{G_{\Bbb{D}}},\Delta _{G_{\Bbb{D}}}\right) =U_{\sigma _0}$ 
as in \emph{\S\ref{TOR-GL}\textsf{(i)}}, we obtain
\[
N_{G_{\Bbb{D}}}={\Bbb{Z}}^d+\Bbb{Z\ }\frak{N}_{\Bbb{D}}, \quad
\left({\Bbb{Z}}^d=\sum_{i=1}^d \Bbb{Z} e_i \right),
\]
where $\left\{ e_1,\ldots ,e_d\right\} $ are the unit vectors of $\left(
N_{G_{\Bbb{D}}}\right) _{\Bbb{R}}\cong \Bbb{R}^d$, and
\[
\frak{N}_{\Bbb{D}}:=\dbigcup\limits_{J\in \frak{D}}\ \left\{ \frac 
1{\mathbf{w}\left( J'\right) }\ \left( e_{\nu _J}-e_{\nu _{J^{\prime
}}}\right) \ \left|
\begin{array}{l}
\text{\emph{for all }}J'\in \frak{D} \\
\text{\emph{with \ }}J'\sqsubset J \\
\text{\emph{and }}\nu _J<\nu _{J'}
\end{array}
\right. \right\} \ .
\]
\emph{(iii) }Furthermore,
$
\#\left( \frak{Y}_{\Bbb{D}}\right) =\#\left( \frak{N}_{\Bbb{D}}\right) =d-1\
.
$
\end{lemma}

\noindent \textit{Proof. }(i) Obviously, each element of the generating
system given in the definition \ref{DEFGD} of~$G_{\Bbb{D}}$ can be written
as product of elements of~$\frak{Y}_{\Bbb{D}}$, and any group generated by a
proper subset of $\frak{Y}_{\Bbb{D}}$ is a proper subgroup of 
$G_{\Bbb{D}}$.\smallskip \newline
(ii) This follows from the usual representation of the group elements of 
$\frak{Y}_{\Bbb{D}}$ by lattice points in~$N_{G_{\Bbb{D}}}$.\smallskip
\newline
(iii) By induction hypothesis, $\#\left( \frak{Y}_{\Bbb{D}_i}\right)
=\#\left( \frak{N}_{\Bbb{D}_i}\right) $ equals $\xi _i-\nu _i$, for all $i$,
$1\leq i\leq p$. Hence,
\begin{eqnarray*}
\#\left( \frak{Y}_{\Bbb{D}}\right) &=&\#\left( \frak{N}_{\Bbb{D}}\right)
=\sum_{i=1}^p\ \left( \xi _i-\nu _i\right) +\left( p-1\right) =\sum_{i=1}^{p-1}\
\xi _i-\sum_{i=2}^p\ \nu _i+(d-1)+\left( p-2\right) \smallskip \\
&=&\sum_{i=1}^{p-1}\ (\nu _{i+1}-1) -\sum_{i=2}^p\  \nu _{i} +
(d-1)+(p-1)=d-1,
\end{eqnarray*}
and the proof is completed. $_{\Box }$
\medskip
\par\noindent
In the following we denote by $\left\{ n_2^{\left( \Bbb{D}\right) },n_3^{\left( 
\Bbb{D}\right) },\ldots ,n_d^{\left( \Bbb{D}\right) }\right\} $ an enumeration of
the $d-1$ elements of $\frak{N}_{\Bbb{D}}$ such that 
\begin{equation}
n_{\nu _i}^{\left( \Bbb{D}\right) }=\frac 1{k_{1,d}}\ \left( e_1-e_{\nu
_i}\right) \ \ \text{ for } 2\leq i\leq p.~\label{NIS}
\end{equation}
Furthermore, we assume that we have chosen an enumeration of the elements $\left\{ 
n_{\nu _i+1}^{\left( \Bbb{D}_i\right)
},\ldots ,\allowbreak n_{\xi _i}^{\left( \Bbb{D}_i\right) }\right\} $ of
$\frak{N}_{\Bbb{D}_i}$ such that 
\[
n_{\nu _i+1}^{\left( \Bbb{D}\right) }=\frac 1{k_{1,d}}\ n_{\nu _i+1}^{\left(
\Bbb{D}_i\right) },\ \ n_{\nu _i+2}^{\left( \Bbb{D}\right) }=\frac
1{k_{1,d}}\ n_{\nu _i+2}^{\left( \Bbb{D}_i\right) },\ \ldots ,\ \ n_{\xi
_i}^{\left( \Bbb{D}\right) }=\frac 1{k_{1,d}}\ n_{\xi _i}^{\left( \Bbb{D}_i\right) 
},\ \ \text{for}\  1\leq i\leq p.
\]

\begin{lemma} \hfill\newline
\vspace{-0.4cm}
\begin{enumerate}
\item[{\rm (i)}] $\left\{e_1, n_2^{\left( \Bbb{D}\right) },n_3^{\left( 
\Bbb{D}\right) },\ldots ,n_d^{\left( \Bbb{D}\right) }\right\}$ is a 
$\Bbb{Z}$-basis of the lattice $N_{G_{\Bbb{D}}}$.
\item[{\rm (ii)}] $\left\{ n_2^{\left( \Bbb{D}\right) },n_3^{\left( \Bbb{D}\right) 
},\ldots ,n_d^{\left( \Bbb{D}\right) }\right\} $ forms a $\Bbb{Z}$-basis of
$
N_{G_{\Bbb{D}}}\cap \emph{lin}\left( \left\{ e_1-e_2,e_1-e_3,\ldots
,e_1-e_d\right\} \right) \ .
$
\end{enumerate}
\label{BASIS}
\end{lemma}

\noindent \textit{Proof. } (i)\, By definition of $N_{G_{\Bbb{D}}}$ it suffices 
to show that 
$\Bbb{Z}^d\subset \Bbb{Z}\{ e_1,n_2^{\left( \Bbb{D}\right) },n_3^{\left( 
\Bbb{D}\right) },\ldots
,n_d^{\left( \Bbb{D}\right) }\}$. 
By induction hypothesis, we may assume 
\[
\Bbb{Z\ }e_{\nu _i}+\Bbb{Z\ }e_{\nu _i+1}+\cdots +\Bbb{Z\ }e_{\xi _i}\subset
\Bbb{Z}\left\{e_{\nu_i}, n_{\nu _i+1}^{\left( \Bbb{D}_i\right)
},\ldots ,\allowbreak n_{\xi _i}^{\left( \Bbb{D}_i\right) }\right\}, 
\ \ \text{for}\  1\leq i\leq p.
\]
Since $k_{1,d}\in \Bbb{Z}$ and $e_{\nu_i}=-k_{1,d}\,n_{\nu_i}^{(\Bbb{D})}+e_1$ we 
also have 
$$
\Bbb{Z\ }e_{\nu _i}+\Bbb{Z\ }e_{\nu _i+1}+\cdots +\Bbb{Z\ }e_{\xi _i}\subset 
\Bbb{Z}\left\{ e_1,n_{\nu_i}^{\left( \Bbb{D}\right) },n_{\nu_{i+1}}^{\left( 
\Bbb{D}\right) },\ldots
,n_{\xi_i}^{\left( \Bbb{D}\right) }\right\}
$$
 and thus 
\[
\bigcup_{i=2}^p\left\{ e_1,n_{\nu _i}^{\left( \Bbb{D}\right) },n_{\nu
_i+1}^{\left( \Bbb{D}\right) },\ldots ,n_{\xi _i}^{\left( \Bbb{D}\right)
}\right\}\ \ =\ \ 
\left\{ e_1,n_2^{\left( \Bbb{D}\right) },
n_3^{\left( \Bbb{D}\right) },\ldots ,n_d^{\left( \Bbb{D}\right)}\right\}
\]
is indeed a $\Bbb{Z}$-basis of $N_{G_{\Bbb{D}}}$. 

\par\noindent
(ii)\, By (i), it suffices to prove that each vector $n_i ^{\left(
\Bbb{D}\right) }$, $2\leq i \leq d$, belongs to the subspace 
lin$( \{e_1-e_2,\ldots ,\allowbreak e_1-e_d\}) $. 
By (\ref{NIS}) this is true for
all $n_{\nu _i}^{\left( \Bbb{D}\right) }$,$\ 2\leq i\leq p$. By induction
hypothesis, we may assume that
\[
n_{\nu _i+j}^{\left( \Bbb{D}_i\right) }\in \text{lin}\left( \left\{ e_{\nu
_i}-e_{\nu _i+1},\ldots ,e_{\nu _i}-e_{\xi _i}\right\} \right),\quad 1\leq i\leq 
p,\,\, 1\leq j\leq \xi _i-\nu _i,
\]
and hence that
\[
n_{\nu _i+j}^{\left( \Bbb{D}\right) }\in \text{lin}\left( \left\{ e_{\nu
_i}-e_{\nu _i+1},\ldots ,e_{\nu _i}-e_{\xi _i}\right\} \right)
\]
too. Since $e_{\nu _i}-e_{\nu _i+j}=\left( e_1-e_{\nu _i+j}\right) -\left(
e_1-e_{\nu _i}\right) $, we are done. $_{\Box }\bigskip $

\noindent $\bullet $ \textit{Third step. }Since the lattice $N_{G_{\Bbb{D}}}$
is ``skew'' and $d$-dimensional, it is not so convenient to work directly
with it. For this reason, we shall perform two affine 
transformations (cf.\ Lemma~\ref{AFITR}), such that $N_{G_{\Bbb{D}}}\cap \text{\rm 
aff}(\frak{s}_{_{G_{\Bbb{D}}}})=N_{G_{\Bbb{D}}}\cap \text{\rm 
aff}(\{e_1,\dots,e_d\})$ becomes the standard lattice $\sum_{i=2}^d 
\Bbb{Z}e_i\cong\Bbb{Z}^{d-1}\subset \Bbb{R}^{d}$. 
\medskip \newline
Define first $\Phi _1^{\Bbb{D}}: \left(N_{G_{\Bbb{D}}}\right)_{\Bbb{R}}\rightarrow 
\left(
N_{G_{\Bbb{D}}}\right) _{\Bbb{R}}$ by 
\[
\Phi _1^{\Bbb{D}}\left( \mathbf{x}\right) :=e_1-\mathbf{x\ .}
\]
Obviously,
\[
\Phi _1^{\Bbb{D}}\left( \frak{s}_{_{G_{\Bbb{D}}}}\right) =\text{conv}\left(
\left\{ \mathbf{0},e_1-e_2,e_1-e_3,\ldots ,e_1-e_d\right\} \right) ,\ \ \ \
\]
and $\Phi _1^{\Bbb{D}}( N_{G_{\Bbb{D}}})=N_{G_{\Bbb{D}}}$.
After that define the affine (linear) transformation  (cf.~Lemma~\ref{BASIS}) 
$\Phi _2^{\Bbb{D}}: \left(N_{G_{\Bbb{D}}}\right)_{\Bbb{R}}  \rightarrow 
\left(N_{G_{\Bbb{D}}}\right)_{\Bbb{R}}$
by setting $\Phi_2^{\Bbb{D}}(e_1)=e_1$ and 
\[
\Phi _2^{\Bbb{D}}\left( \frac 1{\mathbf{w}\left( J'\right) }\
\left( e_{\nu _J}-e_{\nu _{J'}}\right) \right) :=e_{\nu _{J^{\prime
}}},\text{ \ for all elements\ \ }\frac 1{\mathbf{w}\left( J^{\prime
}\right) }\ \left( e_{\nu _J}-e_{\nu _{J'}}\right) \ \text{ of \ 
}\frak{N}_{\Bbb{D}}\ .
\]
Now let $\mathbf{\Phi}^{\Bbb{D}}:\left(N_{G_{\Bbb{D}}}\right)_{\Bbb{R}}  
\rightarrow \left(N_{G_{\Bbb{D}}}\right)_{\Bbb{R}}$ be the composition 
$\mathbf{\Phi}^{\Bbb{D}}:={\Phi}^{\Bbb{D}}_2\circ{\Phi}^{\Bbb{D}}_1$ and let 
$$
\widetilde{\frak{s}_{_{G_{\Bbb{D}}}}}:=\mathbf{\Phi}^{\Bbb{D}}(\frak{s}_{_{G_{\Bbb
{D}}}}) \quad \text {and }\quad 
\widetilde{N_{G_{\Bbb{D}}}^\frak{s}}:=\mathbf{\Phi }^{\Bbb{D}}\left( 
N_{G_{\Bbb{D}}}\cap \text{\rm aff}( \frak{s}_{_{G_{\Bbb{D}}}}) \right).
$$
Observe that $ \widetilde{N_{G_{\Bbb{D}}}^\frak{s}}=\sum_{i=2}^d 
\Bbb{Z}e_i\cong\Bbb{Z}^{d-1}\subset \Bbb{R}^{d}$. 
To provide a convenient description of the
vertices of $\Phi^{\Bbb{D}}(\frak{s}_{_{G_{\Bbb{D}}}})$ incorporating our 
inductive argumentation, we denote the vertex set
of $\widetilde{\frak{s}_{_{G_{\Bbb{D}_i}}}}=         
\mathbf{\Phi}^{\Bbb{D}_i}(\frak{s}_{_{G_{\Bbb{D}_i}}})\subset 
\sum_{j=\nu_i+1}^{\xi_i}\Bbb{R}e_j$ corresponding to the tree 
$\mathbf{W}_{\Bbb{D}_i}$ by 
\[
\text{vert}\left( \widetilde{\frak{s}_{_{G_{\Bbb{D}_i}}}}\right) =\left\{ 
\mathbf{0},\frak{y}_{\nu _i+1}^{\left( \Bbb{D}_i\right) },\frak{y}_{\nu 
_i+2}^{\left(
\Bbb{D}_i\right) },\ldots ,\frak{y}_{\xi _i}^{\left( \Bbb{D}_i\right)
}\right\} ,
\]
in such a way that
\[
\frak{y}_{\nu _i+j}^{\left( \Bbb{D}_i\right) }:=\mathbf{\Phi }^{\Bbb{D}_i}\left( 
e_{\nu _i+j}\right) =\Phi _2^{\Bbb{D}_i}\left( e_{\nu _i}-e_{\nu
_i+j}\right) ,\quad 1\leq i\leq p,\, 1\leq j\leq \xi _i-\nu _i.
\]
Analogously, we label the vertices in
$\text{\rm vert}( 
\widetilde{\frak{s}_{_{G_{\Bbb{D}}}}})=\{\mathbf{0},\frak{y}_{2}^{\left( \Bbb{D} 
\right) },\frak{y}_{3}^{\left(
\Bbb{D}\right) },\ldots ,\frak{y}_{d}^{\left( \Bbb{D}\right)
}\}$ with 
\[
\frak{y}_{j+1}^{\left( \Bbb{D}\right) }:=\mathbf{\Phi }^{\Bbb{D}}\left(
e_{j+1}\right) =\Phi _2^{\Bbb{D}}\left( e_1-e_{j+1}\right) \ \ 
\text{for}\  1\leq j\leq d-1.
\]
\begin{lemma} Using the notation introduced above, 
\begin{eqnarray}
\frak{y}_{\nu _i}^{\left( \Bbb{D}\right) } & = & k_{1,d}\cdot e_{\nu _i},\,\, 
\quad 2\leq i\leq p,  \label{VERT1} \\
\frak{y}_{\nu _1+j}^{\left( \Bbb{D}\right) }& = &k_{1,d}\cdot \frak{y}_{\nu
_1+j}^{\left( \Bbb{D}_1\right) },\quad\,\, 1\leq j\leq \xi_1-\nu_1,   
\label{VERT2}
\\
\frak{y}_{\nu _i+j}^{\left( \Bbb{D}\right) }& = &k_{1,d}\cdot \left( \frak{y}_{\nu 
_i+j}^{\left( \Bbb{D}_i\right) }+e_{\nu _i}\right),\,\,\quad  
2\leq i\leq p,\, \, 1\leq j\leq \xi_i-\nu_i. \label{VERT3}
\end{eqnarray}
\label{hugo}
\end{lemma}
\noindent \textit{Proof. }(\ref{VERT1}) follows from (\ref{NIS})
because
\[
\frak{y}_{\nu _i}^{\left( \Bbb{D}\right) }=\Phi _2^{\Bbb{D}}\left(
e_1-e_{\nu _i}\right) =k_{1,d}\cdot \Phi _2^{\Bbb{D}}\left( n_{\nu
_i}^{\left( \Bbb{D}\right) }\right) =k_{1,d}\cdot e_{\nu _i}\ .
\]
On the other hand,
\[
\Phi _2^{\Bbb{D}}\left( e_{\nu _i}-e_{\nu _i+j}\right) =k_{1,d}\cdot \Phi
_2^{\Bbb{D}_i}\left( e_{\nu _i}-e_{\nu _i+j}\right) =k_{1,d}\cdot \frak{y}_{\nu 
_i+j}^{\left( \Bbb{D}_i\right) }
\]
and (\ref{VERT2}) is clear by setting $i=1$ ($\nu _1=1$). Finally, to get 
(\ref{VERT3}) 
we write
\begin{eqnarray*}
\Phi _2^{\Bbb{D}}\left( e_1-e_{\nu _i+j}\right) &=&\Phi _2^{\Bbb{D}}\left(
\left( e_1-e_{\nu _i}\right) +\left( e_{\nu _i}-e_{\nu _i+j}\right) \right)
\smallskip \\
&=&\Phi _2^{\Bbb{D}}\left( e_1-e_{\nu _i}\right) +\Phi _2^{\Bbb{D}}\left(
e_{\nu _i}-e_{\nu _i+j}\right) =k_{1,d}\cdot e_{\nu _i}+k_{1,d}\cdot \frak{y}_{\nu 
_i+j}^{\left( \Bbb{D}_i\right) }
\end{eqnarray*}
and we are done. $_{\Box }$

\begin{lemma} 
For $\widetilde{\frak{s}_{G_\Bbb{D}}}$ and $\widetilde{N_{G_{\Bbb{D}}}^\frak{s}}$ 
we have: 
\label{DJO}
\begin{enumerate}
\item[{\rm (i)}] \hfill
\vspace{-0.5cm}
$$
\sum_{j=2}^d\Bbb{Z}e_j=\widetilde{N_{G_{\Bbb{D}}}^\frak{s}}=
\text{\rm aff}_{\Bbb{Z}}\left(
\widetilde{N_{G_{\Bbb{D}_1}}^\frak{s}}\cup 
\left(\widetilde{N_{G_{\Bbb{D}_2}}^\frak{s}}+e_{\nu_2}\right)\cup\dots
\cup
\left(\widetilde{N_{G_{\Bbb{D}_p}}^\frak{s}}+e_{\nu_p}\right)
\right).
$$
\item[{\rm (ii)}] \hfill
\vspace{-0.5cm} 
$$
d-1=\dim\left(\widetilde{\frak{s}_{_{G_{\Bbb{D}}}}}\right)=
\dim\left(\widetilde{\frak{s}_{_{G_{\Bbb{D}_1}}}}\right)+
\sum_{i=2}^p 
\dim\left(\widetilde{\frak{s}_{_{G_{\Bbb{D}_i}}}}+e_{\nu_i}\right)+p-1,
$$
\item[{\rm (iii)}] The simplex $\widetilde{\frak{s}_{_{G_{\Bbb{D}}}}}$  can be 
expressed as a dilation of $p-1$ simplex-joins
$$
\widetilde{\frak{s}_{_{G_{\Bbb{D}}}}}= k_{1,d}\ \left( 
\widetilde{\frak{s}_{_{G_{\Bbb{D}_1}}}}*\left( 
\widetilde{\frak{s}_{_{G_{\Bbb{D}_2}}}}+e_{\nu _2}\right) *\cdots *\left(
\widetilde{\frak{s}_{_{G_{\Bbb{D}_p}}}}+e_{\nu _p}\right) \right) 
\subset \sum_{j=2}^d \Bbb{R} e_j \cong \Bbb{R}^{d-1}.
$$

\end{enumerate}
\end{lemma}

\noindent \textit{Proof. } By construction we have 
$\widetilde{N_{G_{\Bbb{D}_i}}^\frak{s}}=\sum_{j=\nu_i+1}^{\xi_i}\Bbb{Z} e_j$ which 
shows (i). Obviously,   
$\dim(\widetilde{\frak{s}_{_{G_{\Bbb{D}_i}}}})=\xi_i-\nu_i$ and since the 
simplices $\widetilde{\frak{s}_{_{G_{\Bbb{D}_1}}}}$, 
$(\widetilde{\frak{s}_{_{G_{\Bbb{D}_2}}}}+e_{\nu_2}),\dots,(\widetilde{\frak{s}_{_
{G_{\Bbb{D}_p}}}}+e_{\nu_p})$ are contained in complementary subspaces, (ii) and 
(iii)  are  immediate consequences of Lemma~\ref{hugo}. $_{\Box }$

\begin{example}
\emph{If we consider the }$7$\emph{-dimensional absolutely unbreakable,}
\emph{msc-,  quotient c.i.-singularity  corresponding to  the special datum
}$\Bbb{D}$
\emph{whose  Watanabe tree  $\mathbf{W}_{\Bbb{D}}$ is shown in
Figure~8,
then the image $\widetilde{\frak{s}_{\Bbb{D}}}$ of the junior
simplex $\frak{s}_{\Bbb{D}}$ under $\mathbf{\Phi }^{\Bbb{D}}$ is the
Watanabe simplex defined as convex hull of the seven vectors:}
\par\noindent
\parbox{8truecm}
{\begin{eqnarray*}
      &       &     (0,0,0,0,0,0,0)^{\intercal}, \\
      &       &     (0,0,{k_{1,7}},0,0,0,0)^{\intercal}, \\
                       &                                             &
(0,0,0,0,{k_{1,7}},{k_{1,7}k_{4,7}k_{4,5}},0)^{\intercal}, \\
      &       &    
(0,0,0,{k_{1,7}},0,{k_{1,7}k_{4,7}},{k_{1,7}k_{4,7}k_{4,5}})^{\intercal}.
\end{eqnarray*}
}
\parbox{7truecm}
{\begin{eqnarray*}
    &  & (0,{k_{1,7}k_{1,2}},0,0,0,0,0)^{\intercal}, \\
    &  & (0,0,0,{k_{1,7}},0,0,0)^{\intercal}, \\
 &  & (0,0,0,{k_{1,7}},0,{k_{1,7}k_{4,7}},0)^{\intercal},
\\
\hphantom{v_8}                              &                             &
\vphantom{(0,0,0,{k_{1,7}},0,{k_{1,7}k_{4,7}},0)^{\intercal},}
\end{eqnarray*}
}
\end{example}
\noindent

\setlength{\GapDepth}{1.0truecm}
\setlength{\GapWidth}{0.1truecm}
\setlength{\EdgeLabelSep}{0.4truecm}
\begin{center}
{\small
\begin{bundle}{\fbox{$\stackrel{ v_{\{1,2,3,4,5,6,7\}} }{ { 1} } $ }}
        \chunk[$k_{1,7}$]{%
        \begin{bundle}{\fbox{$\stackrel{ v_{\{1,2\}} }{ {k_{1,7}}} $ } }
        \chunk[$k_{1,2}$]{\fbox{$\stackrel{ v_{\{1\}} }{ {k_{1,7}k_{1,2}} } $ }  }
        \chunk[$k_{1,2}$]{\fbox{$\stackrel{           v_{\{2\}}          }{
        {k_{1,7}k_{1,2}} } $ } }
        \end{bundle}%
        }
        \chunk[$ k_{1,7}$]{\fbox{$\stackrel{  v_{\{3\}} }{{k_{1,7}}
        } $ } }
        \chunk[$k_{1,7}$]{%
       \begin{bundle}{\fbox{$      \stackrel{       v_{\{4,5,6,7\}}      }{
       {k_{1,7}} } $ } }
        \chunk[$k_{4,7}$]{%
        \begin{bundle}{\fbox{$\stackrel{           v_{\{4,5\}}           }{
       {k_{1,7}k_{4,7}} } $ } }
        \chunk[$k_{4,5}$]{\fbox{$ \stackrel{      v_{\{4\}} }{      {k_{1,7}
        k_{4,7}k_{4,5}} }
        $ } }
        \chunk[$k_{4,5}$]{\fbox{$\stackrel{      v_{\{5\}} }{      {k_{1,7}
        k_{4,7}k_{4,5}} }
        $ } }
        \end{bundle}%
        }
        \chunk[$k_{4,7}$]{%
        \begin{bundle}{\fbox{$\stackrel           {v_{\{6,7\}}           }{
        {k_{1,7}k_{4,7}} } $ } }
        \chunk[$k_{6,7}$]{\fbox{$ \stackrel{     v_{\{6\}} }{      {k_{1,7}
        k_{4,7}k_{6,7}} }
        $ } }
        \chunk[$k_{6,7}$]{\fbox{$       \stackrel{       v_{\{7\}}       }{
        {k_{1,7}k_{4,7}k_{6,7}} }
        $ } }
        \end{bundle}%
        }
\end{bundle}
        }
\end{bundle}}
\vspace{0.1truecm}
\begin{center} {\bf Figure 8} \end{center}
\end{center}

\noindent $\bullet $ \textit{Fourth step. }By Lemma~\ref{DJO} and by
construction, it is now clear that $\widetilde{\frak{s}_{_{G_{\Bbb{D}}}}}$ is a (non-basic) Watanabe $(d-1)$-simplex w.r.t.~$\widetilde{N_{G_{\Bbb{D}}}^\frak{s}}$ and ${\frak{s}_{_{G_{\Bbb{D}}}}}$  a (non-basic) Watanabe $(d-1)$-simplex w.r.t.~${N_{G_{\Bbb{D}}}^\frak{s}}$ whenever $\mathbf{W}_{\Bbb{D}}$ is a non-trivial Watanabe tree.
 If $\mathbf{W}_{\Bbb{D}}$ is a Watanabe forest,
then we recall the ``superscripts in square brackets'' from \ref{SPLIT}
(iii): Taking its connected components to be $\mathbf{W}_{\Bbb{D}^{\left[
i\right] }}$, $1\leq i\leq \kappa $, the group $G_{\Bbb{D}}$ splits into the
direct product of $G_{\Bbb{D}^{\left[ i\right] }}$'s, and by arguments
similar to those of the previous steps we can show that 
$\frak{s}_{_{G_{\Bbb{D}}}}$ is again a transformed join of 
$\frak{s}_{_{G_{\Bbb{D}^{\left[
1\right] }}}}$ with the junior simplices $\frak{s}_{_{G_{\Bbb{D}^{\left[
i\right] }}}}$, $2\leq i\leq \kappa $, corresponding to the other
components, after translating by unit vectors. (We do not need an extra
dilation in this case). Hence, the one direction of the
Reduction Theorem \ref{RED} 
is completely proved. The proof of the converse statement is based on
the ``backtracking method'', and on the fact that any toric affine variety
being associated to a rational s.c.p.~cone  supported by  a (non-basic) simplex, 
is the
underlying space of an abelian quotient singularity; it is therefore
omitted.

\begin{example}
\emph{Up to isomorphism, the  only four-dimensional abelian, msc-, quotient
c.i.-singularity  which  ``breaks''  is  that  corresponding  to the datum}
$\Bbb{D}$ \emph{and Watanabe forest} $\mathbf{W}_{\Bbb{D}}$ 
\emph{of Figure~9.  (In fact,  it splits  into two  Hirzebruch-Jung 
singularities of types}
$A_{k_{1,2}-1}$\emph{\ and} $A_{k_{3,4}-1}$\emph{, respectively).}
\vspace{-.4truecm}
\end{example}

%
\setlength{\GapDepth}{0.8truecm}
\setlength{\GapWidth}{0.2truecm}
\setlength{\EdgeLabelSep}{0.4truecm}
\begin{center}
\hfill\hbox{
\begin{bundle}{\fbox{$ v_{\{1,2\}}, { 1} $ }}
        \chunk[$ k_{1,2}$]{\fbox{$ v_{\{1\}}, {k_{1,2}} $ } }
        \chunk[$ k_{1,2}$]{\fbox{$ v_{\{2\}}, {k_{1,2}} $ } }
\end{bundle}
\hskip 2truecm
\begin{bundle}{\fbox{$ v_{\{3,4\}}, { 1} $ }}
        \chunk[$ k_{3,4}$]{\fbox{$ v_{\{3\}}, {k_{3,4}} $ } }
        \chunk[$ k_{3,4}$]{\fbox{$ v_{\{4\}}, {k_{3,4}} $ } }
\end{bundle}
}\hfill
\vspace{0.1truecm}
\begin{center} {\bf Figure 9} \end{center}
\end{center}

\noindent The junior tetrahedron $\frak{s}_{_{G_{\Bbb{D}}}}$ (which is the
join of two $1$-dimensional Watanabe simplices), is drawn in Figure 10, for the
values of free parameters $k_{1,2}=3$\emph{, } $k_{3,4}=5$.

\begin{center}
\unitlength 1.3truecm
\begin{picture}(2,1.8)
\thicklines
\jput(0,0){$\bullet$} \put(-.3,0){$\mathbf{e_1}$}
\jput(2,0){$\bullet$} \put(2.3,0){$\mathbf{e_2}$}
\jput(1,1.732){$\bullet$} \put(1.2,1.77){$\mathbf{e_4}$}
\jput(1,0.7){$\bullet$} \put(0.72,0.75){$\mathbf{e_3}$}
\put(0.666,0){$\bullet$}
\put(1.333,0){$\bullet$}
\put(1,0.9064){$\bullet$}
\put(1,1.1128){$\bullet$}
\put(1,1.3192){$\bullet$}
\put(1,1.5256){$\bullet$}
\put(0.05,0.05){\drawline(0,0)(1,1.732)(2,0)(0,0)}
\thinlines
\put(0.05,0.05){\drawline(0,0)(1,0.7)(2,0)(0,0)(1,0.7)(1,1.732)}
\end{picture}
\vspace{-0.2truecm}
\begin{center} {\bf Figure 10} \end{center}
\end{center}

\section{Final step: Basic, coherent triangulations of the junior simplex}

\noindent In this section we give the proof of our Main Theorem \ref{Main}.
To obtain the desired basic, coherent triangulations of the junior simplices
$\frak{s}_G$ of all abelian quotient c.i.-singularities, it suffices (by \ref
{AFITR}, \ref{RED}) to construct b.c.b.-triangulations for all
Watanabe simplices. (For simplicity's sake, in the proofs of \ref{bass} and \ref{Watti}, we shall assume that our reference lattice is the standard lattice $\Bbb{Z}^d$.)

\begin{proposition}
\label{bass}All dilations $\lambda \,\mathbf{s}$
for integral $\lambda \geq 2$, 
of a basic lattice $d$-simplex $\mathbf{s}$ in $\Bbb{R}^d$,
have b.c.b.-triangulations.
\end{proposition}

\noindent \textit{Proof. }Up to an affine integral transformation, we may assume 
that $\mathbf{s}$ equals $\mathbf{s}_d$ (w.r.t.\ $\Bbb{Z}^d$), where
\[
\mathbf{s}_d:=\text{conv}\left( \left\{ \mathbf{0},e_1,e_1+e_2,\ldots
,e_1+\cdots +e_d\right\} \right) =\left\{ \mathbf{x\in \,}\Bbb{R}^d\,\left|
\,0\leq x_d\leq x_{d-1}\leq \ldots \leq x_1\leq 1\right. \right\} \ .
\]
The affine hyperplane arrangement $\Bbb{H}_d$ (of type 
$\widetilde{\mathcal{A}}_d$) consisting of the union of hyperplanes
\begin{eqnarray*}
& &      \mathcal{H}_i\left(      k\right)     =\left\{     \mathbf{x\in
\,}\Bbb{R}^d\,\left|   \,x_i=k\right.   \right\},
\ \ \ \text{for}\ 1\leq   i\leq   d,\ \ k\in
\Bbb{Z}, \\
\text{\rm and} & &\\
& &  \mathcal{H}_{i,j}\left( k\right) =\left\{ \mathbf{x\in \,}\Bbb{R}^d\,\left| 
\,x_i-x_j=k\right. \right\}
\ \ \ \text{for}\ 1\leq i<j\leq d,\ \ k\in \Bbb{Z},
\end{eqnarray*}
is infinite, but only the hyperplanes
\begin{eqnarray*}
& &      \mathcal{H}_i\left(      k\right)     =\left\{     \mathbf{x\in
\,}\Bbb{R}^d\,\left|   \,x_i=k\right.   \right\},
\ \ \ \text{for}\ 1\leq   i\leq   d,\ \ 1\le k\le \lambda-1 \\
\text{\rm and} & &\\
& &  \mathcal{H}_{i,j}\left( k\right) =\left\{ \mathbf{x\in \,}\Bbb{R}^d\,\left| 
\,x_i-x_j=k\right. \right\}
\ \ \ \text{for}\ 1\leq i<j\leq d,\ \ 1\le k\le \lambda-1
\end{eqnarray*}
intersect the interior of $\lambda \,\mathbf{s}_d$. The hyperplanes 
$\mathcal{H}_i\left( k\right) ,1\leq i\leq d,k\in \Bbb{Z}$, subdivide $\Bbb{R}^d$ 
into unit cubes. Given such a cube
\[
C\left( \mathbf{\mu }\right) :=\left[ 0,1\right] ^d+\mathbf{\mu }\text{ \ \
\ \ (where }\mathbf{\mu }=\left( \mu _1,\ldots ,\mu _d\right) \in 
\Bbb{Z}^d\text{)}
\]
and fixing $i,j$ with $1\leq i<j\leq d$, we see that there is only one of
the hyperplanes $\mathcal{H}_{i,j}\left( k\right) $ intersecting the
interior of $C\left( \mathbf{\mu }\right) $, namely $\mathcal{H}_{i,j}\left(
\mu _i-\mu _j\right) $. The hyperplanes
\[
\left\{ \mathcal{H}_{i,j}\left( \mu _i-\mu _j\right) \ \left| \ 1\leq
i<j\leq d\right. \right\}
\]
provide the ``usual'' basic triangulation of 
$C\left( \mathbf{\mu }\right) $ into $d!$ basic subsimplices of the form
\[
s\left( \mathbf{\mu ,}\theta \right) =\mathbf{\mu }+\text{conv}\left(
\left\{ \mathbf{0},e_{\theta \left( 1\right) },e_{\theta \left( 1\right)
+\theta \left( 2\right) },\ldots ,e_{\theta \left( 1\right) +\theta \left(
2\right) +\cdots +\theta \left( d\right) }\right\} \right) ,\ \text{for 
all permutations\ }\theta \in \frak{S}_d.
\]
Thus, $\Bbb{H}_d$ defines a basic triangulation $\frak{T}_{\left( d\right) }$
of the entire space $\Bbb{R}^d$ (cf.\ \cite[Ch.~III]{KKMS}). This
triangulation is also coherent because
\[
\Bbb{R}^d\ni \mathbf{x}
\longmapsto \overline{\psi
}\left( \mathbf{x}\right) =-
\sum_{0\leq i<j\leq d}
\left\{ 
\sum_{0\le k\le x_j-x_i}    H(x_j{-}x_i{-}k) \ +\ 
\sum_{x_j-x_i\le k\le 0}    H(k{-}x_j{+}x_i)
\right\}
 \in \Bbb{R}
\]
with $x_0:=0$, using the Heaviside function $H:\Bbb{R}\rightarrow\Bbb{R}$ 
\[
H(x):=\left\{
\begin{array}{cl}
x & \text{\rm if } x\geq 0, \\
0 & \text{\rm otherwise,}
\end{array}
\right.
\]
is a strictly upper convex function for this triangulation.
Finally, it is balanced
via the colouring function $\varphi _d$:
\[
\Bbb{Z}^d
\ni \mathbf{\mu }=\left( \mu _1,\ldots ,\mu _d\right) \ \stackrel{\varphi 
_d}{\mathbf{\mapsto
}}\mathbf{\ }\left( \mu _1+\cdots +\mu _d\right) \ \text{mod }\left(
d+1\right)  \in \left\{ 0,1,\ldots ,d\right\}
\]
which $\left( d+1\right) $-colours all the facets $s\left( \mathbf{\mu ,}\theta 
\right) $ of $\frak{T}_{\left( d\right) }$. 

Since the bounding
hyperplanes of $\lambda \,\mathbf{s}_d$ are contained in $\Bbb{H}_d$, the
restriction $\frak{T}_{\left[ d;\lambda \right] }:=\frak{T}_{\left( d\right)
}\left| _{\lambda \,\mathbf{s}_d}\right. $ is a b.c.b.-triangulation of $\lambda 
\,\mathbf{s}_d$.
 $_{\Box }$

\begin{example}
\emph{Figure 11 provides the ``nice'' basic triangulation }$\frak{T}_{\left[
2;4\right] }$ \emph{of }$4\,\mathbf{s}_2$ \emph{inherited from the affine
hyperplane arrangement} $\Bbb{H}_2$.
\end{example}

\begin{center}
\unitlength 0.7truecm
\begin{picture}(5,5)
\multiput(0,0)(0,1){5}{\multiput(0,0)(1,0){5}{\circle*{0.15}}}
\thicklines
\put(0,0){\line(1,0){4}}
\put(0,0){\line(1,1){4}}
\put(4,0){\line(0,1){4}}
\thinlines
\dottedline{0.1}(-1,0)(5,0)
\dottedline{0.1}(-1,1)(5,1)
\dottedline{0.1}(-1,2)(5,2)
\dottedline{0.1}(-1,3)(5,3)
\dottedline{0.1}(-1,4)(5,4)
\dottedline{0.1}(-1.5,3.5)(0,5)
\dottedline{0.1}(-1.5,2.5)(1,5)
\dottedline{0.1}(-1.5,1.5)(2,5)
\dottedline{0.1}(-1.5,0.5)(3,5)
\dottedline{0.1}(-1.5,-0.5)(4,5)
\dottedline{0.1}(-1,-1)(5,5)
\dottedline{0.1}(0,-1)(5.5,4.5)
\dottedline{0.1}(1,-1)(5.5,3.5)
\dottedline{0.1}(2,-1)(5.5,2.5)
\dottedline{0.1}(3,-1)(5.5,1.5)
\dottedline{0.1}(0,-1)(0,5)
\dottedline{0.1}(1,-1)(1,5)
\dottedline{0.1}(2,-1)(2,5)
\dottedline{0.1}(3,-1)(3,5)
\dottedline{0.1}(4,-1)(4,5)
\end{picture}
\vspace{0.8truecm}
\begin{center} {\bf Figure 11} \end{center}
\end{center}

\begin{corollary}
\label{COR-H}Theorem \emph{\ref{Main} }is true for all $\left( d;k\right) 
$-hypersurface-singularities.
\end{corollary}

\noindent \textit{Proof. }In this case the junior simplex $\frak{s}_G$
equals (up to an affine integral transformation) $k\,\mathbf{s}$, where
\[
\mathbf{s}=\text{conv}\left( \left\{ \frac 1k\ e_1,\ldots ,\frac 1k\
e_d\right\} \right)
\]
w.r.t.\ $N_{G\left( d;k\right) }$. But $\mathbf{s}$ is obviously basic. $_{\Box }$

\begin{remark}
\emph{Corollary \ref{COR-H} generalizes Roan's result in 
\cite[\S\ 5]{Roan1}. (Roan proved the existence of crepant
desingularizations of }$\left( d;2\right) $\emph{-hypersurfaces by
successively blowing up the singular parts of the branching locus of
the corresponding double covering of }$\Bbb{C}^d$; \emph{see 
Remark~\ref{REMH}(i)).}
\end{remark}

\noindent 
By Theorem~\ref{RED}, and what we have already explained
in \S\ \ref{TOR-CR-R}, our Main Theorem \ref{Main} is now a
consequence of the following:

\begin{theorem}\label{Watti}
Every Watanabe simplex has a b.c.b.-triangulation.
\end{theorem}

\noindent \textit{Proof. }Trivially, every Watanabe $0$-simplex has a
b.c.b.-triangulation. If $\mathcal{T}_{1}$ (resp.\ $\mathcal{T}_{2}$)
denotes a b.c.b.-triangulation of a Watanabe $d_1$-simplex $\mathbf{s}_1$
(resp.\ of a Watanabe $d_2$-simplex $\mathbf{s}_2$), then 
$\mathcal{T}_{1}*\mathcal{T}_{2}$ is a b.c.b.-triangulation of 
$\mathbf{s}_1*\mathbf{s}_2$
by Theorem~\ref{JOINS}. Hence, it suffices to show that for a given Watanabe
simplex $\mathbf{s}$, being equipped with a b.c.b.-triangulation, say 
$\mathcal{T}$, having a colouring function 
$\varphi : \mathbf{s} \cap \Bbb{Z}^d\rightarrow \left\{ 0,1,\ldots ,d\right\} $ 
and a
strictly upper convex function $\psi : \mathbf{s}  \rightarrow
\Bbb{R}$, the $\lambda $-times dilation of $\mathbf{s}$ possesses itself a
b.c.b.-triangulation. Every facet $F$ of $\mathcal{T}$ gives rise to a
\textit{unique} affine integral transformation $\Phi _F:\mathbf{s}_d\rightarrow F$ 
($\mathbf{s}_d$ as before in the proof of \ref{bass}) which
respects the colourings relatively to $\mathbf{s}_d$, i.e., which sends every
vertex of $\mathbf{s}_d$ to the vertex of $F$ that has same colour. $\Phi _F$
maps $\lambda \,\mathbf{s}_d$ onto $\lambda \,F$. Thus, the image
\[
\widetilde{\mathcal{T}_{\lambda \,F}}:=\Phi _F\left( \frak{T}_{\left[
d;\lambda \right] }\right) =\lambda \,\frak{T}_{\left( d\right) }\left|
_{\lambda \,F}\right.
\]
constitutes a triangulation of $\lambda \,F$ (which is a
b.c.b.-triangulation by \ref{dilations}, \ref{AFITR}, and \ref{bass}). The
triangulations $\left\{ \widetilde{\mathcal{T}_{\lambda \,F}}\ \left| \ F\text{ 
facets of }\mathcal{T}\right. \right\} $ fit together to give a basic
triangulation $\widehat{\mathcal{T}}$ of $\lambda \,\mathbf{s}$. 
 Since $\widehat{\psi }:\left| \widehat{\mathcal{T}}\right| 
 \rightarrow \Bbb{R}$
 defined by
 \[
 \lambda \,F\ni \mathbf{x}\mapsto \widehat{\psi}\left( \mathbf{x}\right)
 :=\psi \left( \frac 1\lambda \ \mathbf{x}\right) +\varepsilon \cdot
 \overline{\psi }\left( \left( \Phi_F \right) ^{-1}\left( \mathbf{x}\right)
 \right) \in \Bbb{R}
 \]
(for $\varepsilon>0$ sufficiently small and fixed for {\it all} facets $F$ of $\mathcal{T}$, and for $\overline{\psi}$ as in \ref{bass}), 
 is strictly upper convex,  
the
triangulation $\widehat{\mathcal{T}}$ is also coherent by Patching Lemma \ref
{PATCH}. Finally, $\widehat{\mathcal{T}}$ is balanced because 
$\widehat{\varphi }:| \widehat{\mathcal{T}}| \rightarrow \Bbb{R}$ 
defined by
\[
\lambda \,F\ni \mathbf{x}\mapsto \widehat{\varphi}\left( \mathbf{x}\right) 
:=\varphi _d\left( \left( \Phi _F\right) ^{-1}\left( \mathbf{x}\right) \right) \in 
\Bbb{R}
\]
for \textit{all} facets $F$ of $\mathcal{T}$ (with $\varphi _d$ as in \ref{bass}), is a colouring map. $_{\Box }$

\begin{example}
\emph{Let }$\mathcal{T}$ \emph{denote the unique b.c.b.-triangulation of the
}$2$\emph{-dimensional Watanabe simplex }
\[
\mathbf{s}=\emph{conv}\left( \left\{ \left(
\begin{array}{r}
0 \\
0
\end{array}
\right) ,\left(
\begin{array}{r}
-1 \\
1
\end{array}
\right) ,\left(
\begin{array}{c}
-1 \\
-1
\end{array}
\right) \right\} \right) \emph{\ \ \ (w.r.t.\ }\Bbb{Z}^2\emph{)}
\]
{\rm with the facets}
$F_1=\text{\rm conv}(\{{-1\choose 1},{-1\choose 0},{0\choose 0}\})$ {\rm and}  
$F_2=\text{\rm conv}(\{{-1\choose-1},{-1\choose 0},{0\choose 0}\})$.
\emph{\ Figure 12 illustrates how one ``glues together'' }
\[
\widetilde{\mathcal{T}_{3F_i}}=\Phi _{F_i}\left( \frak{T}_{\left[ 2;3\right]
}\right) =3\,\frak{T}_{\left( 2\right) }\left| _{3\,F_i}\right. \emph{,\ \ 
}i=1,2\emph{,}
\]
\emph{\ to obtain the b.c.b.-triangulation }$\widehat{\mathcal{T}}$. \emph{The 
affine integral transformations} $\Phi _{F_i}$
\emph{are given by}
\[
\Phi _{F_1}\left( \mathbf{x}\right) =\left(
\begin{array}{rr}
 0 & 1 \\
 1 & 0
\end{array}
\right) \mathbf{x}+\left(
\begin{array}{r}
-1 \\
-1
\end{array}
\right) ,\ \ \ \Phi _{F_2}\left( \mathbf{x}\right) =\left(
\begin{array}{rr}
 0 & 1 \\
-1 & 0
\end{array}
\right) \mathbf{x}+\left(
\begin{array}{r}
-1 \\
 1
\end{array}
\right) ,\ \ \mathrm{for}\ \mathbf{x}\in \Bbb{R}^2.
\]
\end{example}

\[
\unitlength1cm
\begin{picture}(0,0)
\put(3.8,1.1){$\mathbf{s}_2$}
\put(-.5,0.3){the basic simplex $\mathbf{s}_2$}
\put(-.6,9.2){the triangulation $\mathcal{T}$}
\put(11,9){the entire}
\put(11,8.5){gluing triangulation $\widehat{\mathcal{T}}$}
\put(11,0.1){the ``nice'' triangulation}
\put(11,-.4){$\frak{T}_{[2;3]}$ of $3\mathbf{s}_2$}
\put(2.5,7.3){$F_1$}
\put(2.5,7.8){$F_2$}
\put(4,3){$\Phi_{F_1}$}
\put(0.5,3){$\Phi_{F_2}$}
\put(12,6){$\Phi_{F_1}$}
\put(6.8,3){$\Phi_{F_2}$}
\put(6,-1){\textbf{Figure 12}}
\end{picture}
\def\epsfsize#1#2{.75#1}
\epsffile{pastry2.eps}
\]

\vspace{1cm}

\section{On the computation of cohomology group dimensions}

\noindent Let $\Bbb{C}^d/G=U_{\sigma _0}=X\left( N_G,\Delta _G\right) $ be
again the underlying space of an abelian quotient c.i.-singularity and
\begin{equation}
f_{\mathcal{T}}:X\left( N_G,\widehat{\Delta _G}\left( \mathcal{T}\right)
\right) \rightarrow X\left( N_G,\Delta _G\right)  \label{RESOL}
\end{equation}
any partial crepant desingularization induced by a lattice triangulation 
$\mathcal{T}$ of the junior simplex $\frak{s}_G$. It is easy to verify that
the central fiber $\mathbf{F}_{\mathcal{T}}=\left( f_{\mathcal{T}}\right)
^{-1}\left( \left[ \mathbf{0}\right] \right) $ of $f_{\mathcal{T}}$ is a
strong deformation retract of the overlying space $X\left( N_G,\widehat{\Delta 
_G}\left( \mathcal{T}\right) \right) $. Theorem \ref{Main} guarantees
the existence of \textit{at least one }b.c.b.-triangulation $\mathcal{T}$ of
$\frak{s}_G$, so that (\ref{RESOL}) is a projective, crepant, full
desingularization. But even if we let $\mathcal{T}$ go through the entire
class of \textit{all} possible basic triangulations of $\frak{s}_G$, and use
\cite[Thm.~5.4]{BD}, \cite[Cor.~1.5]{Ito-Reid}, we obtain the one-to-one
McKay-type correspondence
\[
\left\{ \text{elements of \ }G\text{ }\right\} \stackrel{1:1}{\longleftrightarrow 
}\left\{
\begin{array}{c}
\text{a basis of \ }H^{*}\left( \mathbf{F}_{\mathcal{T}},\Bbb{Q}\right)
\text{ consisting } \\
\text{ of classes of algebraic cycles}
\end{array}
\right\}
\]
and, in particular,
\vspace{-4truemm}
\[
\left\{ \text{elements of \ }G\text{ \ of age\emph{\ }}1\right\} 
\stackrel{1:1}{\longleftrightarrow }\left\{
\begin{array}{c}
\text{exceptional prime} \\
\text{divisors w.r.t.\emph{\ }}f_{\mathcal{T}}
\end{array}
\right\} \stackrel{1:1}{\longleftrightarrow }\left\{
\begin{array}{c}
\text{a basis of \ }H^2\left( \mathbf{F}_{\mathcal{T}},\Bbb{Q}\right)\\
\text{ consisting of classes of } \\
\text{algebraic cycles}
\end{array}
\right\}
\]
where the elements of age $1$ are those having
lattice-point-representatives lying on $\frak{s}_G$. In fact, only the even
cohomology groups of $\mathbf{F}_{\mathcal{T}}$ are non-trivial. To compute
their dimensions we need some concepts from enumerative
combinatorics.\medskip \newline
For a lattice $d$-polytope $P\subset N_{\Bbb{R}}$ w.r.t.\ an $N\cong \Bbb{Z}^d $, 
and $\kappa $ a positive integer, let

\[
\mathbf{Ehr}_N\left( P,\kappa \right) =\mathbf{a}_0\left( P\right) 
+\mathbf{a}_1\left( P\right) \ \kappa +\cdots +\mathbf{a}_{d-1}\left( P\right) \
\kappa ^{d-1}+\mathbf{a}_d\left( P\right) \ \kappa ^d\ \in \ \Bbb{Q}\left[
\kappa \right]
\]
denote the \textit{Ehrhart polynomial }of $P$ (w.r.t.\ $N$), where 
$\mathbf{Ehr}_N\left( P,\kappa \right) :=\#\,\left( \kappa \,P\cap N\right) $, and
\[
\mathbf{Ehr}_N\left( P;\frak{t}\right) :=1+\sum_{\kappa =1}^\infty \ 
\mathbf{Ehr}_N\left( P,\kappa \right) \ \frak{t}^\kappa \in \Bbb{Q}_{\,}\left[
\!\left[ \frak{t}\right] \!\right]
\]
\textit{\ }the corresponding \textit{Ehrhart series}. Writing 
$\mathbf{Ehr}_N\left( P;\frak{t}\right) $ as
\[
\mathbf{Ehr}_N\left( P;\frak{t}\right) =\frac{\ \mathbf{\delta }_0\left(
P\right) +\ \mathbf{\delta }_1\left( P\right) \ \frak{t}+\cdots +\ \mathbf{\delta 
}_{d-1}\left( P\right) \ \frak{t}^{d-1}+\ \mathbf{\delta }_d\left(
P\right) \ \frak{t}^d}{\left( 1-\frak{t}\right) ^{d+1}}
\]
we get the so-called $\mathbf{\delta }$-\textit{vector} $\mathbf{\delta }\left( 
P\right) =\left( \mathbf{\delta }_0\left( P\right) ,\mathbf{\delta }_1\left( 
P\right) ,\ldots ,\mathbf{\delta }_{d-1}\left( P\right) ,\mathbf{\delta }_d\left( 
P\right) \right) $ of $P.$

\begin{definition}
\emph{For any integer} $d\geq 0$ \emph{we introduce the} \textit{transfer} 
$\mathbf{a}$-$\mathbf{\delta }$-\textit{matrix }$\mathcal{M}_d\in $ 
\emph{GL}$\left( d+1,\Bbb{Q}\right) $ \emph{(depending only on }$d$\emph{) to be
defined as}
\[
\mathcal{M}_d:=\left( \frak{R}_{i,j}\right) _{0\leq i,j\leq d}\ \ \ \text{\emph{\ 
with}}\ \ \ \ \frak{R}_{i,j}:=\dfrac 1{d!}\ \left\{
\dsum\limits_{p=i}^d\ \QDATOPD[ ] {d}{p}\ \dbinom pi\ \left( d-j\right)
^{p-i}\right\}
\]
\emph{where} $\QDATOPD[ ] {d}{p}$ \emph{denotes the Stirling number (of the
first kind) of} $d$ \emph{over} $p$.
\end{definition}


\noindent The following lemma can be proved easily.

\begin{lemma}
For a lattice $d$-polytope $P\subset N_{\Bbb{R}}$ w.r.t.\ an $N\cong \Bbb{Z}^d $, 
we have
\[
\left( \mathbf{a}_0\left( P\right) ,\mathbf{a}_1\left( P\right) ,\ldots 
,\mathbf{a}_{d-1}\left( P\right) ,\mathbf{a}_d\left( P\right) \right) =\left(
\mathbf{\delta }_0\left( P\right) ,\mathbf{\delta }_1\left( P\right) ,\ldots
,\mathbf{\delta }_{d-1}\left( P\right) ,\mathbf{\delta }_d\left( P\right)
\right) \cdot \left( \mathcal{M}_d\right) ^{\intercal}.
\]
\end{lemma}

\begin{theorem}
For any integer $d\geq 2$, and for any crepant \emph{(}full\emph{) 
}desingularization \emph{(\ref{RESOL}) }of an abelian quotient c.i.-space 
$\Bbb{C}^d/G$, the non-trivial cohomology dimensions of $\mathbf{F}_{\mathcal{T}}$ 
\emph{(}or of $X\left( N_G,\widehat{\Delta _G}\left( \mathcal{T}\right)
\right) $, for any basic $\mathcal{T}$\emph{)}, can be determined
inductively \emph{(}via the coefficients of Ehrhart polynomials\emph{)} as
follows \emph{:\smallskip }\newline
\emph{(i) }If $\frak{s}_G$ 
is the join of two Watanabe simplices $\mathbf{s}_1$, $\mathbf{s}_2$, of 
dimensions $d_1$, $d_2$, \emph{(}$d_1+d_2=d-2\emph{)}$,
then

$
\text{\emph{dim}}_{\Bbb{Q}}H^{2i}\left( \mathbf{F}_{\mathcal{T}},\Bbb{Q}\right) =
$
\begin{equation}
=\sum\Sb 0\leq p,q\leq d-1  \\ \text{\emph{with \thinspace }}p+q=i  \endSb 
\,\left( \left( \left[ \left( \mathcal{M}_{d-1}\right) ^{-1}\right]
_{p+1}\cdot \left(
\begin{array}{c}
\mathbf{a}_0\left( \mathbf{s}_1\right) \smallskip \\
\,\mathbf{a}_1\left( \mathbf{s}_1\right) \smallskip \\
\vdots \smallskip \\
\mathbf{a}_{d-1}\left( \mathbf{s}_1\right)
\end{array}
\right) \right) \times \left( \left[ \left( \mathcal{M}_{d-1}\right)
^{-1}\right] _{q+1}\cdot \left(
\begin{array}{c}
\mathbf{a}_0\left( \mathbf{s}_2\right) \smallskip \\
\mathbf{a}_1\left( \mathbf{s}_2\right) \smallskip \\
\vdots \smallskip \\
\mathbf{a}_{d-1}\left( \mathbf{s}_2\right)
\end{array}
\right) \right) \right)  \label{KOHJ}
\end{equation}
\ \smallskip \ \newline
\emph{(ii) }If $\frak{s}_G$ \emph{(}up to an affine integral transformation\emph{) 
}is the dilation $\lambda \,\mathbf{s}$, $\lambda \geq 2$, of a
Watanabe simplex $\mathbf{s}$, then
\begin{equation}
\text{\emph{dim}}_{\Bbb{Q}}H^{2i}\left( \mathbf{F}_{\mathcal{T}},\Bbb{Q}\right) 
=\left[ \left( \mathcal{M}_{d-1}\right) ^{-1}\right] _{i+1}\cdot
\left(
\begin{array}{c}
\mathbf{a}_0\left( \mathbf{s}\right) \smallskip \\
\lambda \,\mathbf{a}_1\left( \mathbf{s}\right) \smallskip \\
\vdots \smallskip \\
\lambda ^{d-1}\,\mathbf{a}_{d-1}\left( \mathbf{s}\right)
\end{array}
\right)  \label{KOH-D}
\end{equation}
for all $i$, $0\leq i\leq d-1$, where for all $\rho $, $1\leq \rho \leq d$, 
$\left[ \left( \mathcal{M}_{d-1}\right) ^{-1}\right] _\rho $ denotes the $\rho 
$-th row-vector of~$\left( \mathcal{M}_{d-1}\right) ^{-1}$.
\end{theorem}

\noindent \textit{Proof. }By \cite[Thm.~4.4]{BD}, dim$_{\Bbb{Q}}H^{2i}\left(
\mathbf{F}_{\mathcal{T}},\Bbb{Q}\right) $ equals the $i$-th component of the
$\mathbf{\delta }$-vector of $\frak{s}_G$.
\smallskip \newline
That the computation
can be done inductively follows from the converse statement in Theorem~\ref{RED}
and Theorem~\ref{Main}. In fact, for any $G$ being conjugate to a $G_{\Bbb{D}}$, 
the cohomology dimensions can be read off from the free parameters of 
$\mathbf{W}_{\Bbb{D}}$.
\smallskip \newline
(i) Since $\mathcal{T}$ is basic, the $\mathbf{\delta }$-vector of $\frak{s}_G$ is 
equal to the $h$-vector of  $\mathcal{T}$ (see Stanley \cite[2.5]{St-CP}). Hence, (\ref
{KOHJ}) follows from the formula which provides the $h$-vector of the join
of two simplicial complexes.\smallskip \newline
(ii) (\ref{KOH-D}) is obvious because $\mathbf{a}_j\left( \lambda 
\,\mathbf{s}\right) =\lambda ^j\,\mathbf{a}_j\left( \mathbf{s}\right) $ for all 
$j$, $0\leq j\leq d-1$. $_{\Box }$


\begin{corollary}
Let $G$ be conjugate to $G\left( d;k\right) $ \emph{(}within \emph{SL}$\left( 
d,\Bbb{C}\right) $\emph{)}. Then the non-trivial cohomology
dimensions of any crepant resolution \emph{(\ref{RESOL}) }of the $\left(
d;k\right) $-hypersurfaces equals

\begin{equation}
\text{\emph{dim}}_{\Bbb{Q}}H^{2i}\left( \mathbf{F}_{\mathcal{T}},\Bbb{Q}\right) 
=\left[ \left( \mathcal{M}_{d-1}\right) ^{-1}\right] _{i+1}\cdot
\left(
\begin{array}{c}
1 \\
\vdots \\
\stackunder{\text{\emph{entry of the column at} }\left( j+1\right) 
\text{\emph{-position} \emph{(}}0\leq j\leq d-1\text{\emph{)} 
}}{\underbrace{\dfrac{k^j}{\left( d-1\right) !}\ \ \left\{ 
\dsum\limits_{p=j}^{d-1}\
\QDATOPD[ ] {d-1}{p}\ \dbinom pj\ \left( d-1\right) ^{p-j}\right\} }} \\
\vdots \smallskip \\
\dfrac{k^{d-1}}{\left( d-1\right) !}
\end{array}
\right)  \label{KOH-H}
\end{equation}
for $0\leq i\leq d-1$. In particular, the Euler-Poincar\'{e}
characteristic $\chi \left( \mathbf{F}_{\mathcal{T}}\right) $ of 
$\mathbf{F}_{\mathcal{T}}$ equals
\begin{equation}
\chi \left( \mathbf{F}_{\mathcal{T}}\right) =\chi \left( X\left( N_{G\left(
d;k\right) },\widehat{\Delta _{G\left( d;k\right) }}\left( \mathcal{T}\right) 
\right) \right) =\left| G\left( d;k\right) \right| =k^{d-1}
\label{KOH-EP}
\end{equation}
\end{corollary}

\noindent \textit{Proof. }In this case (up to an affine integral
transformation) $\frak{s}_G$ equals $k\,\mathbf{s}$, where
\[
\mathbf{s}=\text{conv}\left( \left\{ \frac 1k\ e_1,\ldots ,\frac 1k\
e_d\right\} \right)
\]
w.r.t.\ $N_{G\left( d;k\right) }$. Since
\[
\mathbf{Ehr}_{N_{G\left( d;k\right) }}\left( \mathbf{s},\kappa \right) 
=\dbinom{\kappa +d-1}\kappa =\dbinom{\kappa +d-1}{d-1}=\frac 1{\left(
d-1\right) !}\ \sum_{p=0}^{d-1}\ \QDATOPD[ ] {d-1}{p}\ \left( \kappa
+d-1\right) ^p\,,\
\]
or, alternatively, since
\[
\mathbf{\delta }_0\left( \mathbf{s}\right) =1\text{ \ \ \ and \ \ \ 
}\mathbf{\delta }_1\left( \mathbf{s}\right) =\cdots =\mathbf{\delta }_{d-1}\left(
\mathbf{s}\right) =0,
\]
we get
\[
\mathbf{a}_j\left( \mathbf{s}\right) =\left\{
\begin{array}{c}
\text{ }\left( j+1\right) \text{-entry of the} \\
\text{first column of} \\
\text{the matrix }\mathcal{M}_{d-1}
\end{array}
\right\} =\frac 1{\left( d-1\right) !}\ \ \left\{ \dsum\limits_{p=j}^{d-1}\
\QDATOPD[ ] {d-1}{p}\ \dbinom pj\ \left( d-1\right) ^{p-j}\right\}
\]
for all $j$, $0\leq j\leq d-1$, and (\ref{KOH-H}) follows from (\ref{KOH-D}). 
Finally, formula (\ref{KOH-EP}) follows from the equality $\chi \left(
\mathbf{F}_{\mathcal{T}}\right) =\left( d-1\right) !\,\mathbf{a}_{d-1}\left(
\mathbf{s}\right) $ . $_{\Box }$

\section{Comments and open problems.}

\noindent (i) All partial, crepant, projective desingularizations (\ref
{RESOL}) of $\left( \Bbb{C}^d/G,\left[ \mathbf{0}\right] \right) =\left(
U_{\sigma _0},\text{orb}\left( \sigma _0\right) \right) $ can be studied by
means of the \textit{secondary polytope} of the junior simplex $\frak{s}_G$,
whose vertices parametrize all coherent triangulations $\mathcal{T}$ of 
$\frak{s}_G$. Passing from one vertex of this polytope to another, we perform
a finite series of flops. It should be mentioned that for $d\geq 4$, it is
possible to start from a vertex corresponding to a basic $\mathcal{T}$, and
arrive at another maximal, but non-basic triangulation $\mathcal{T}^{\prime
} $. In the language of ``toric MMP'' (see Reid \cite{Reid2}), our
b.c.b.-triangulations lead to \textit{smooth }minimal models. (Note that all
morphisms $f_{\mathcal{T}}$ can be written, by \cite[0.2-0.3]{Reid2}, as
compositions of finite sequences of more elementary toric
contraction-morphisms).\medskip

\noindent (ii) For any $d\geq 3$, the general hyperplane section of the
above $\left( \Bbb{C}^d/G,\left[ \mathbf{0}\right] \right) $ through $\left[ 
\mathbf{0}\right] $ is either a rational or an elliptic Gorenstein
singularity (see \cite[\S~3.10]{Reid3}). Already for $d=3$, both
possibilities occur. For example, the general hyperplane section of the $\left( 
3;k\right) $-hypersurfaces is a Du Val singularity (of type $D_4$) for 
$k=2$, and an elliptic Gorenstein surface singularity whenever $k\geq 3$. It might 
be interesting to investigate the class of Gorenstein
elliptic singularities obtained by this procedure by exploiting the
inductive character of Watanabe's classification, essentially via the forests 
$\mathbf{W}_{\Bbb{D}}.$ (What is the relationship between these
general-hyperplane-section singularities and the free parameters of the
starting-point singularities for $d\geq 4\,$?)\medskip

\noindent (iii) Could Theorem~\ref{Main} be generalized for the underlying
spaces of Gorenstein, toric, non-quotient, c.i.-singularities? More
precisely, what would be the geometric analogue of joins and dilations
describing the structure of lattice polytopes which support the
Gorenstein cones in this case?\medskip

\noindent (iv) Theorem \ref{Main} has various applications to global
geometrical constructions. For instance, every ``well-stratified''
Calabi-Yau variety which is locally a complete intersection, and has at most
abelian quotient singularities, possesses global, crepant, full resolutions
in all dimensions. (Nevertheless, to check the projectivity of these
globally desingularizing morphisms one needs to apply the Nakai-Moishezon
criterion, and this is only possible if one has some extra information
available about the intrinsic geometry of the varieties being under
consideration).\medskip

\noindent (v) A special class of such Calabi-Yau varieties, which is of
particular interest, is that of compactified hypersurfaces $\overline{Z_f}$ 
$\subset Y_P$ being embedded in a toric variety $Y_P$, associated to a 
\textit{reflexive}, \textit{simple} lattice polytope $P$. Assuming that 
$\overline{Z_f}$ is $P$\textit{-regular }(in Batyrev's sense \cite{B}), and
that the polar lattice polytope $P^{\,\star }$ of $P$ has only
Watanabe simplices as faces of codimension $\geq 2$, there exist always
global, crepant, full desingularizations $\widehat{Z_f}\rightarrow \overline{Z_f}$ 
of $\overline{Z_f}$. One method to construct at least one of them is
to triangulate the faces of $P^{\,\star }$ as in the present paper and then
to join the single interior point of $P^{\,\star }$ with them. For example,
the mirror-partner $\overline{Z_g}$ of a marginally deformed
Fermat-hypersurface $\overline{Z_f}$ (or of any smooth hypersurface) of
degree $d$ in $\Bbb{P}_{\Bbb{C}}^{d-1}$ has global, crepant, full 
desingularizations
in \textit{all} dimensions. (In this particular case, one can, in addition,
easily construct a globally projective desingularizing morphism). Hence, the
so-called ``string-theoretic'' Hodge numbers of this $\overline{Z_g}$ (cf.\ 
\cite{B} \cite{BD}) 
are nothing but the \textit{usual }Hodge-numbers of an always
existing \textit{smooth, projective} $\widehat{Z_g}$; in particular, the
corresponding monomial-divisor mirror-map provides the usual dualism between
the polynomial first-order deformations of $\overline{Z_f}$ within 
$H^{d-3,1}\left( \overline{Z_f}\right) $ on the one hand, and the toric part
of $H^{1,1}\left( \widehat{Z_g}\right) $, on the other. (So there are
concrete exceptional divisors, and there is no need here to work in the
category of singular spaces). This motivates the formulation of a purely
combinatorial problem: to \textit{classify} all reflexive, simplicial
polytopes, at least in dimension $5$, having only Watanabe simplices as
faces of codimension at least $2$.
\bigskip \medskip

\noindent \textit{Acknowledgements. }
The first author would like express his thanks to  F. Hirzebruch
and M.~Reid for useful discussions and remarks, to DFG for
the support by an one-year-fellowship, 
and to the Mathematics Institute of Bonn University for hospitality.\\ 
The second and third author  wish to acknowledge the long-term support 
by a DFG Gerhard-Hess Forschungsf\"{o}rderpreis. Since October 1996 the second 
author is supported by the DFG Leibniz-Preis of M.~Gr{\"o}tschel.


\end{document}